\documentclass[a4paper,12pt]{article}

\usepackage{subeqnarray,amsfonts,graphicx,enumerate}
\usepackage[squaren]{SIunits} \usepackage{color,xspace}
\usepackage{subfig} \newcommand{\eps}{\varepsilon}
\usepackage{amsmath,amssymb} \usepackage{paralist}
\usepackage{graphics} \usepackage{epsfig}
\usepackage[colorlinks=true]{hyperref}
\newcommand{\beq}{\begin{equation}}
\newcommand{\eeq}{\end{equation}}
\newcommand{\eq}[1]{(\ref{#1})}

\DeclareMathOperator{\minusone}{-}

    \newcommand{\R}{\field{R}}
    \newcommand{\field}[1]{\mathbb{#1}}

    \newcommand{\ie}{\textit{i.e.}, }
    \newcommand{\eg}{\textit{e.g.}, }

    \newcommand{\RefSec}[1]{Sec.~\ref{#1}}
    
    \newcommand{\InRefFig}[1]{Figure~\ref{#1}}
    
    \newcommand{\RefFig}[1]{Fig.~\ref{#1}}
    
    \newcommand{\RefEq}[1]{Eq.~(\ref{#1})}
    
    \newcommand{\Ref}[1]{(\ref{#1})}

    \newlength{\FigureWidthSmall}
    \newlength{\FigureTwoFrames}
    \newlength{\ArraySkip}
    \newlength{\ArraySmallSkip}
    \newlength{\ArrayBigSkip}
    \newlength{\ResColW}
    \newlength{\ResColWS}
    \newlength{\TabHeadSkip}
    \newlength{\TabEndSkip}
\title{Codimension-two homoclinic bifurcations underlying spike adding in the Hindmarsh-Rose burster}
\author{Daniele Linaro${}^1$, Alan Champneys${}^2$, Mathieu Desroches${}^2$, Marco Storace${}^1$\\
\small{${}^1$Biophysical and Electronic Engineering Department, University of Genoa, Genova, Italy}\\
\small{${}^2$Department of Engineering Mathematics, University of Bristol, Bristol, UK}}
\date{\today}

\begin{document}

\maketitle
\begin{abstract}

The well-studied Hindmarsh-Rose model of neural action potential is revisited from the point of
view of global bifurcation analysis.  This slow-fast system of three paremeterised differential
equations is arguably the simplest reduction of Hodgkin-Huxley models capable of exhibiting all
qualitatively important distinct kinds of spiking and bursting behaviour.  First, keeping the
singular perturbation parameter fixed, a comprehensive two-parameter bifurcation diagram is
computed by brute force. Of particular concern is the parameter regime where lobe-shaped regions of
irregular bursting undergo a transition to stripe-shaped regions of periodic bursting. The boundary
of each stripe represents a fold bifurcation that causes a smooth spike-adding
transition where the number of spikes in each burst is increased by one.
Next, numerical continuation studies reveal that the global structure is
organised by various curves of homoclinic bifurcations.

In particular the lobe to stripe transition is organised by a sequence of codimension-two orbit-
and inclination-flip points that occur along {\em each} homoclinic branch. Each branch undergoes a
sharp turning point and hence approximately has a double-cover of the same curve in parameter
space. The sharp turn is explained in terms of the interaction between a two-dimensional unstable
manifold and a one-dimensional slow manifold in the singular limit. Finally, a new local analysis
is undertaken using approximate Poincar\'{e} maps to show that the turning point on each homoclinic
branch in turn induces an inclination flip that gives birth to the fold curve that organises the
spike-adding transition.  Implications of this mechanism for explaining spike-adding behaviour in
other excitable systems are discussed.
\end{abstract}

\section{Introduction}

The Hindmarsh-Rose (HR) model \cite{Hindmarsh:1984} is one of the most widely studied parameterised
three-dimensional systems of ordinary differential equations (ODEs) that arises as a reduction of
the conductance-based Hodgkin-Huxley model for neural spiking \cite{Hodgkin:1952}.  Its success
comes from both its simplicity --- just three ODEs with polynomial nonlinearity, and only a few key
parameters --- and its ability to qualitatively capture the three main dynamical behaviours
displayed by real neurons, namely quiescence, tonic spiking and bursting (see \RefFig{fig:behav}).
Moreover, transitions between these behaviours can be easily described in terms of the
bio-physically motivated parameters. Even reduced-order models like the HR equations can have
direct physiological meaning and so can be used to match or indeed predict detailed \textit{in vivo}
recordings; for instance, in \cite{deLange:2008} the authors use the HR model -- after
appropriately rescaling the state variable $x$, the parameter $I$ and time -- to fit the activity
of both pyramidal cells and neocortical interneurons. Nevertheless, a key argument for their use is
that they can point to generic understanding of which kinds of interventions or perturbations are
likely to lead to certain kinds of transition. These understandings can then be used to help guide
parameter searches for more in-depth computational models which can only be investigated by direct
numerical simulation (DNS). In turn, these simulations can help guide experimental or clinical
control strategies or protocols.

Many papers have investigated the bifurcations that occur in the HR model upon variation of one or
more of its parameters; see
\cite{Belykh:2000,Gonzalez:2003,Gonzalez:2007,Innocenti:2009,Innocenti:2007,Shilnikov:2008,Storace:2008:Chaos,Terman:1991,Terman:1992,Wang:1993}.
These studies have typically focussed on particular transitions: from periodic to irregular
(chaotic) spiking-bursting dynamics, from tonic spiking to bursting and on the two possible kinds
of bursting (square-wave and pseudo-plateau). For perhaps the most comprehensive bifurcation
analysis to date the reader is referred to the work of Shilnikov \& Kolomiets
\cite{Shilnikov:2008}.

In this paper we shall be concerned with understanding the complete bifurcation scenario that
underlies the {\em smooth} transition from tonic spiking to bursting, paying particular attention
to an observed sequence of spike-adding transitions.  This form of period-adding behaviour would
cause a variation of the average number of spikes within a burst, which behaviour is believed to
have important physiological implications \cite{Osinga:2010}.  The key point of the paper is to
show that codimension-one homoclinic bifurcations and their degeneracies are crucial to
understanding how such transitions are organised in parameter space. The methodology we shall adopt
will be a combination of brute-force methods (augmenting the preliminary results in
\cite{Storace:2008:Chaos}), slow-fast arguments, numerical continuation (using AUTO07P
\cite{Auto07PMan}), and geometric analysis using approximate Poincar\'e maps.

Brute-force methods involve classification of stable asymptotic behaviour computed using DNS over a
wide range of parameter values. When trying to understand the mechanisms behind observed
transitions in stable behaviour though, such methods can only describe bifurcation scenarios
roughly. In particular, they often fail to uncover coexisting attractors, and they do not compute
the unstable invariant sets that are often crucial in isolating and analysing both local and global
bifurcations.

From a more analytic point of view, the HR model can be decomposed into a reduced two-dimensional
``fast'' ODE-system with an additional slow variable. Such slow-fast arguments, see
\cite{Bertram95,Golubitsky01, Izhikevich:2000,Rinzel85,Shilnikov:2008}, can provide much generic
information about the original model and tend to work best close to the singular limit of infinite
time-scale separation. However, most physical systems operate away from the singular limit, and the
mutual interactions between slow and fast dynamics are typically very subtle and give rise to
further bifurcations in the complete system that occur ``beyond all orders'' in the singular limit.

Numerical continuation analysis~\cite{krauskopf07} is typically quite robust and through its use of
boundary-value problems to solve for recurrent trajectories does not suffer from the same problems
as DNS in the singular limit. Nevertheless, as we shall see, problems can still arise in the
presence of ``canard-like'' phenomena~\cite{Desroches,Guckenheimer}. In this case, a mix of
numerical results and geometrical analysis can prove pivotal.

The rest of this paper is organised as follows. The next section introduces the Hindmarsh-Rose
model and its dynamics, in particular highlighting its slow-fast structure. Section \ref{sec:3}
then presents numerical bifurcation analysis computed both by brute force and through numerical
continuation. Particular attention is focussed on an infinite family of bifurcation curves of
homoclinic orbits that connect an equilibrium on the unstable part of the critical manifold to itself. Various codimension-two homoclinic bifurcation points are detected and local bifurcation curves
arising from them are computed. It is found that the key fold bifurcations underlying spike-adding
transitions originate from sharp turning point along the homoclinic bifurcation curves, caused by
the interaction between a 1D slow manifold and a 2D unstable manifold of an equilibrium point. Owing
to the sharpness of the folding of the unstable manifold, numerical continuation is found to be
inconclusive. Section \ref{sec:4} therefore presents a new geometric analysis of such a situation
and shows that there has to be an additional codimension-two homoclinic bifurcation of inclination
flip type very close to this point of interaction. Finally section \ref{sec:5} draws conclusions,
suggests avenues for further work, and points to some wider implications of the results.

\section{The Hindmarsh-Rose model}
\label{sec:2}

\begin{figure}[!t]
\centering
\includegraphics[width=10cm]{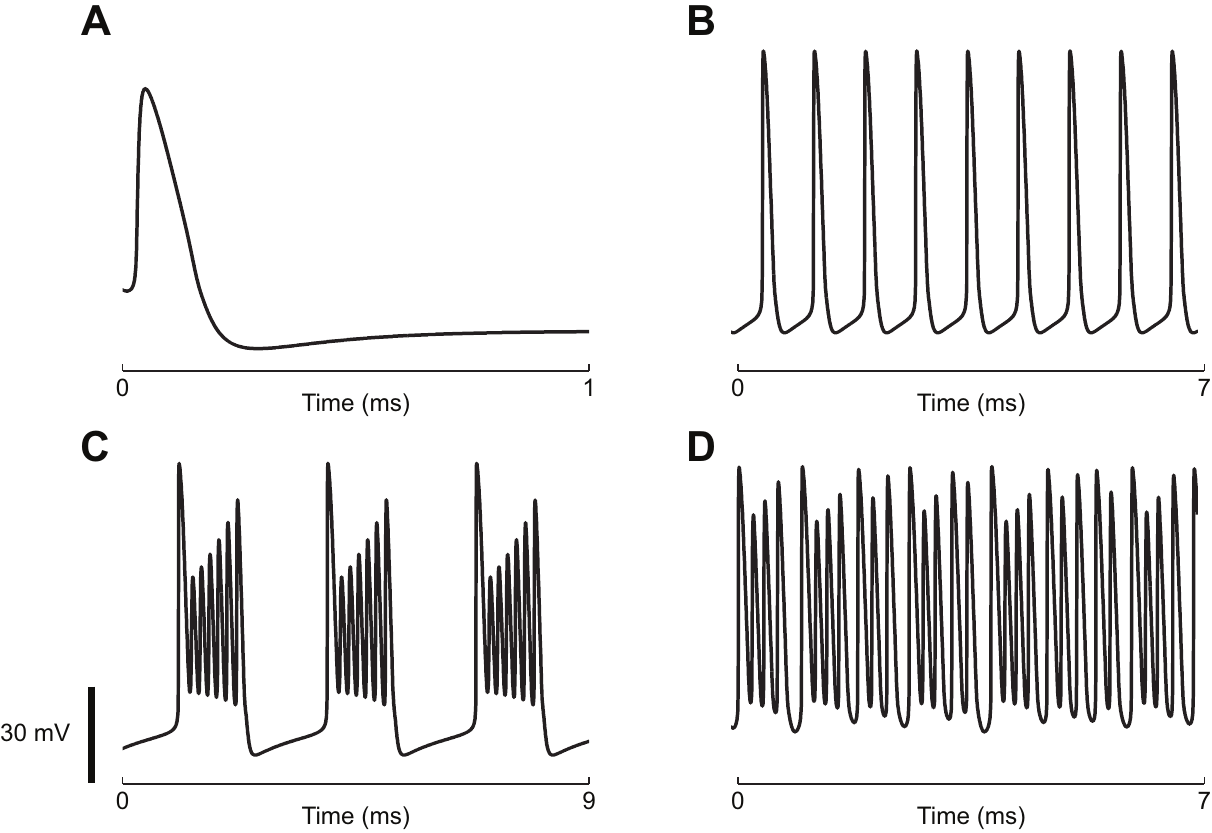}
\caption[Typical neuronal behaviours]{Typical membrane potentials computed
using a model of a leech heart interneuron
showing, under varying input parameters,
(A) quiescence, (B) spiking, (C) bursting and (D) irregular spiking.
The model is as in \cite{Shilnikov:2005}, with
parameters $C=0.5$, $E_K=-0.07$, $E_{Na}=0.045$, $g_{Na}=200$,
$g_1=8$, $E_1=-0.046$, $\tau_{Na}=0.0405022$, $V_{K2}^{shift}=-0.0145$
and $g_{K2}=30$ and: (A) $I_{app}=-0.02$,
$\tau_{K2}=0.25$; (B) $I_{app}=0$, $\tau_{K2}=0.2$; (C) $I_{app}=0$,
$\tau_{K2}=1$;  (D), $I_{app}=0.07$, $\tau_{K2}=0.35$.
}
\label{fig:behav}
\end{figure}

Typical spiking neurons occurring across biology can undergo a variety of distinct dynamical
behaviours, according to the values of bio-physical parameters. See \RefFig{fig:behav} for
outputs from a biologically relevant neuron model.  Among the most important behaviours, one may
find \cite{Izhikevich:2000}:
\begin{itemize}
\item \textit{quiescence}, which occurs when the input to the neuron is below a certain
    threshold and the output does not include any action potential firing events (or spikes);
\item \textit{spiking}, in which the output is made up of a regular series of spikes;
\item \textit{bursting}, where the output consists of groups of two or more spikes separated by
    periods of inactivity;
\item \textit{irregular spiking}: the output is made up of an aperiodic series of spikes.
\end{itemize}

Previous studies have shown that the HR model is able to reproduce all these dynamical behaviours
\cite{Belykh:2000,Gonzalez:2003,Gonzalez:2007,Innocenti:2009,Innocenti:2007,Shilnikov:2008,Storace:2008:Chaos,Terman:1991,Terman:1992,Wang:1993}.
Moreover, in these references bifurcation analysis has been carried out in one or two parameters,
unfolding cascades of smooth transitions between stable bursting solutions and continuous spiking
regimes, both regular (periodic) and irregular (chaotic). See Fig.~\ref{fig:HRbifwide} below for a
recapitulation of some of these results in a two-parameter bifurcation diagram.

\subsection{The governing equations}

The phenomenological neuron model proposed by Hindmarsh and Rose (HR)
\cite{Hindmarsh:1982,Hindmarsh:1984} is a single-compartment model that is computationally simple
yet is capable of mimicking the rich variety of firing pattern behaviours exhibited by real
biological neurons \cite{Herz:2006}. It can be described by the following set of ODEs:
\begin{equation} \label{eq:HR}
\left\{
\begin{array}{@{}r@{}c@{}l}
\dot x &=& y -x^3 + b x^2 + I - z \\[0.5em]
\dot y &=& 1 - 5 x^2 - y \\[0.5em]
\dot z &=& \mu \left( s \left( x - x_{rest} \right) - z \right)
\end{array}
\right.
\end{equation}
The model is dimensionless and the variables have only phenomenological interpretations.  The
variables $x$ and $y$ represent the fast charging dynamics (voltage and current respectively)
associated with a single neuron whereas $z$ is a slow variable mirroring the action of slow ionic
channels.  Hence, \eqref{eq:HR} is a slow-fast system with two fast and one slow variable. Its fast
nullcline $M_{eq}:=\{(x,y,z)\in \mathbb{R}^3;\;z=y-x^3+bx^2+I,\;y=1-5x^2\}$ is the so-called
\textit{critical manifold} of the system. The critical manifold $M_{eq}$ is a manifold of
equilibria for the limiting problem obtained by setting $\mu=0$ in \eqref{eq:HR}, and plays a
crucial role in the non-trivial dynamics of the full system; see section~\ref{slowfast} below. The
roles played by the system parameters can be described as follows. The parameter $I$ mimics the
membrane input current for biological neurons, whereas $b$ is an excitability parameter that allows
one to switch between bursting and spiking behaviours and to control the spiking frequency. The
variable $\mu$ controls the time scale of the slow variable $z$, that is, the
efficiency of the slow channels in exchanging ions. In the presence of spiking behaviour, it
affects the inter-spike interval, whereas in the case of bursting it affects the number of
spikes per burst. The phenomenological parameter $s$ governs the degree of adaptation in the
neuron.  A value of $s$ around unity causes spiking behaviour with no spike-frequency accommodation
nor subthreshold adaptation, whereas values around $s=4$ (the value we shall use in this paper)
allow strong accommodation and subthreshold overshoot, and can even allow oscillations. The
parameter $x_{rest}$ sets the resting potential of the system and is usually set to $-1.6$ in the
dimensionless units in which \Ref{eq:HR} is written.

In what follows, unless otherwise stated we shall consider \Ref{eq:HR} at parameter values
\begin{equation}
\label{eq:parvalues}
s=4, \qquad x_{rest}=-1.6, \qquad \mu=0.01
\end{equation}
and allow $I$ and $b$ to be bifurcation parameters.

\subsection{Spike adding and slow-fast analysis}
\label{slowfast} 
\begin{figure}[!b]
\centering
\includegraphics[width=10cm]{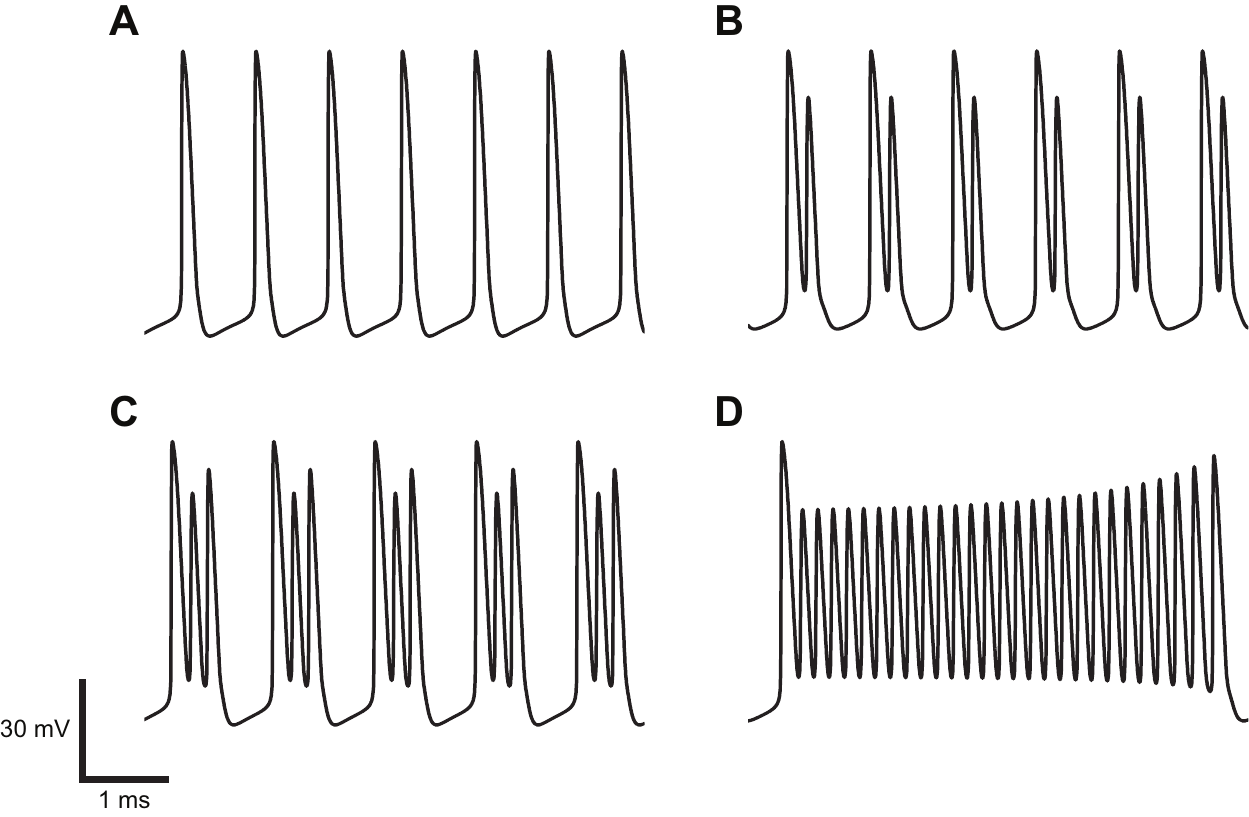}
\caption[Qualitative example of period adding]
{Examples of period adding in the same model as in Fig.~\ref{fig:behav}. By varying a single
parameter, periodic responses can be observed which undergo a transition from (A) isolated
spikes, to (B) two spikes per burst, to (C) three spikes per burst, and so on, up to
(D) 28 spikes per burst.
Model and parameters as in Fig.~\ref{fig:behav} but with
$I_{app}=0$ and $\tau_{K2}=0.2$, $0.3$, $0.5$ and $3.5$
for panels A to D respectively.
}
\label{fig:pa}
\end{figure}

A key feature displayed by many neuronal models is the so-called spike adding
mechanism: this means that it is possible, by changing one parameter, to {\em smoothly} change the
behaviour of the system from spiking to bursting with increasing number of spikes per bursts, as
shown in \RefFig{fig:pa}. In the HR model, changing the parameter $b$ allows one to observe this
transition among periodic responses.  To some extent, also raising the injected current $I$ leads
to this spike adding, but this also typically increases frequency of the bursts. Other combinations
of parameters, \eg $\mu$ and $I$, lead to similar results
\cite{Gonzalez:2007,Rech:2011,Storace:2008:Chaos}.

\begin{figure}[!b]
\centering
\includegraphics[width=10cm]{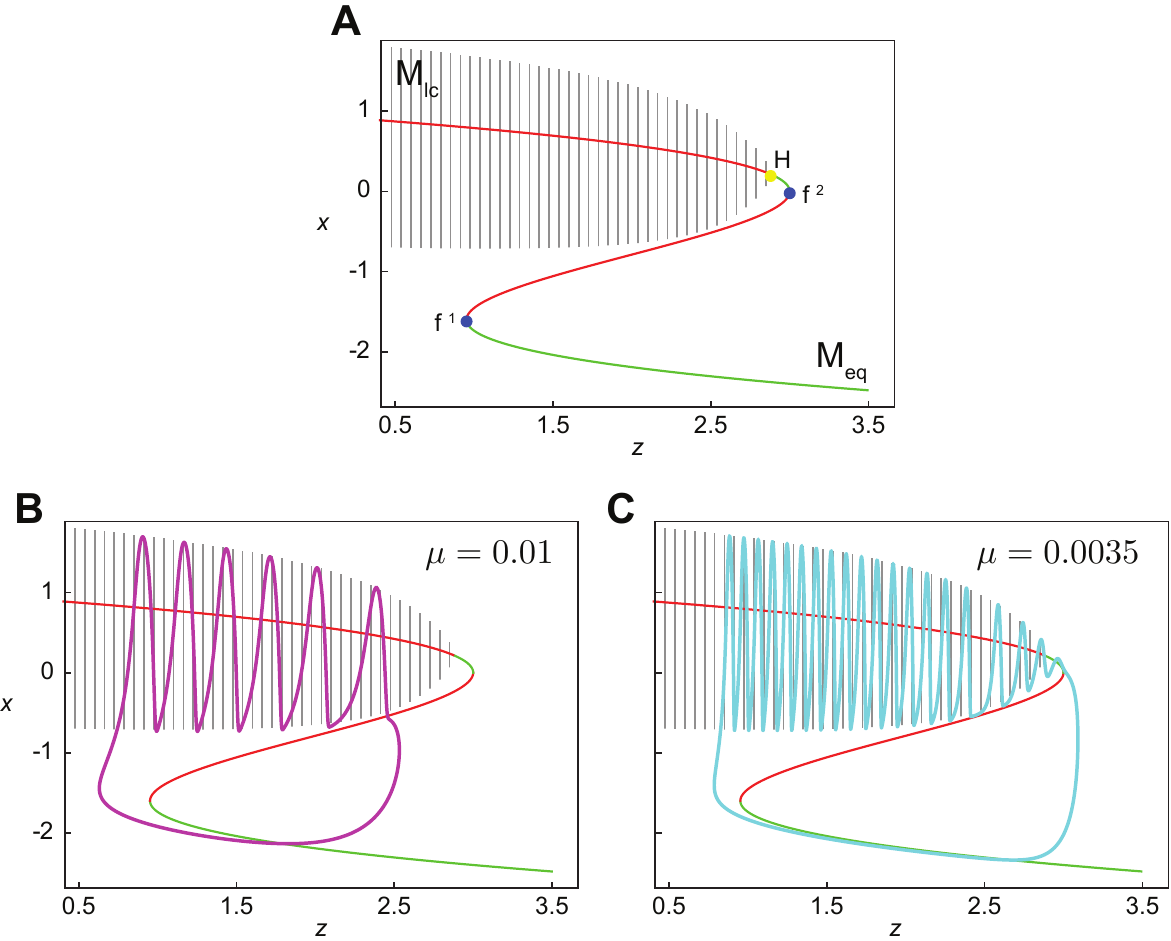}
\caption[Slow-fast analysis of bursting in the HR neuron model]
{Slow-fast decomposition of bursting in the HR neuron model with $b=2.6$ and $I=2$.
Panel A shows the bifurcation
analysis of the fast subsystem obtained by varying the slow variable $z$ considered as a parameter. Panels B and C show
solutions of the full HR model for different values of $\mu$: it is evident that the value of
$\mu$ influences the number of turns of the solution around the manifold of limit cycles, \ie
the number of spikes per burst. See text for more details.}
\label{fig:slowfast}
\end{figure}

One of the most widely used approaches to analyse period-adding in the HR and related models is
so-called {\em slow-fast analysis}
\cite{Desroches11,Guckenheimer,Shilnikov:2005,Shilnikov:2005:PRL,Shilnikov:2008,Terman:1991}. The
technique consists of separating the model into two or more ODE subsets, corresponding to parts of
the neuron that operate on different time scales, for example membrane voltage and fast currents on
the one hand and slow currents and calcium dynamics on the other hand. We shall now provide an
overview of such analysis in system~\eqref{eq:HR}. The {\em fast subsystem} is given by the
equations
\begin{equation} \label{eq:HRfast}
\left\{
\begin{array}{l}
\dot x = y - x^3 + b x^2 + I - z \\
\dot y = 1 - 5 x^2 - y.
\end{array}
\right.
\end{equation}
which contain only two out of the five initial parameters ($b$ and $I$), but in the limit $\mu=0$,
$z$ becomes a constant parameter. Note that the critical manifold $M_{eq}$ of system~\eqref{eq:HR}
corresponds to the set of equilibrium points of the fast subsystem~\eqref{eq:HRfast}. To illustrate
this point, \RefFig{fig:slowfast}A depicts a one-parameter bifurcation analysis of system
\Ref{eq:HRfast} obtained by varying $z$ and keeping $b=2.6$ and $I=2$ fixed. There is a curve of
equilibria (shown in green in the figure) which is stable for large $z$ and undergoes two folds
$f^{1,2}$, becomes stable again before loosing stability at a supercritical Hopf bifurcation $H^-$.
The grey vertical lines are the projections onto the $(z,x)$ plane of the stable limit cycles of
system~\eqref{eq:HRfast}, which organise the {\em bursting} behaviour of the full HR model. It is
interesting to note that for the particular $b$ and $I$ values chosen, at around $z=2.5$ the stable
limit cycle becomes close to the unstable equilibrium that exists between $f^1$ and $f^2$. The
system is therefore close to a homoclinic bifurcation, and indeed we observed the typical
drop-shaped cycles which are characteristic of homoclinic trajectories.

The results of this simple bifurcation analysis for $\mu=0$ constitute the {\em skeleton} of the
solution of the full HR model (for $\mu
>0$). In particular, for sufficiently small $\mu$ geometric singular
perturbation theory \cite{Fenichel:1979} gives, away from the non-hyperbolic points --- here, fold
points --- of the critical manifold $M_{eq}$ (as equilibria of the fast subsystem) the existence of
centre-like (hence, non-unique) one-dimensional perturbed (locally invariant) manifolds
$O(\mu)$-close to the manifold $M_{eq}$ of equilibria. Families of hyperbolic equilibria of the
fast dynamics are referred to as \textit{normally hyperbolic invariant manifolds} for the full
system~\cite{Jones95}. These perturbed manifolds are called \textit{slow (or Fenichel) manifolds}.
Fenichel theory \cite{Fenichel:1971,Fenichel:1979} accounts for their existence and behaviour near
hyperbolic equilibrium points of the fast dynamics only; however, close to non-hyperbolic points,
they can be extended by the flow and their behaviour can be understood using other analytical
techniques such as \textit{geometric
  desingularization} or
\textit{blow-up}~\cite{Krupa:2001}. Furthermore, generalised Fenichel theory~\cite{Fenichel:1971}
also accounts for the persistence of manifolds of fast motion from the singular limit $\mu=0$ to
the $0<\mu\ll 1$ case. In the particular case we are investigating here, the family of limit cycles
of the fast subsystem, parameterised by the (frozen) slow variable $z$, is hyperbolic and, hence,
persists as a two-dimensional manifold $M_{lc}$ close to these limit cycles.

A simple explanation of bursting in the full model is that the solution repeatedly \emph{switches}
between $M_{eq}$ and $M_{lc}$ under the action of the slow variable $z$. This leads to periods of activity, the bursts, in which the solution
evolves close to $M_{lc}$ interspersed with periods of quiescence, in which the solution moves close to
$M_{eq}$. To better clarify this concept, panels B and C of \RefFig{fig:slowfast} show the
projections onto the $(z,x)$-plane of two bursting solutions of the complete system for different
values of $\mu$. For $\mu=0.01$ (panel B), the number of spikes per burst is 6, whereas
$\mu=0.0035$ (panel C) leads to bursts with 20 spikes (panel C). The increased number of spikes is
due to the fact that, smaller $\mu$ causes slower trajectories along the $z$-direction. As a
consequence, the solution stays for a longer period of time close to $M_{lc}$ and, because the {\em fast}
rotation within $M_{lc}$ happens at a $\mu$-independent rate, the number of spikes per burst
increases as $\mu$ decreases. Note that the fast subsystem~\eqref{eq:HRfast} is by definition independent of $\mu$, so one can superimpose the attractor of the full system onto the bifurcation diagram of~\eqref{eq:HRfast} for any value of $\mu$. As we shall see, depending on the values of all the other parameters,
these hybrid bursting trajectories can be either periodic, quasi-periodic or chaotic. In parameter
regions where periodic responses are found, the mechanism by which new spikes are added is complex
and involves interaction of the bursting region with the fold $f^1$ of the critical manifold.  The main
purpose of this paper then is to provide a bifurcation theoretic explanation of period-adding
process. In particular we shall seek an explanation that works for general $\mu$-values without
reference to slow-fast analysis.

\section{Numerical bifurcation analysis}
\label{sec:3}

\begin{figure}[!t]
\centering
\includegraphics[width=10cm]{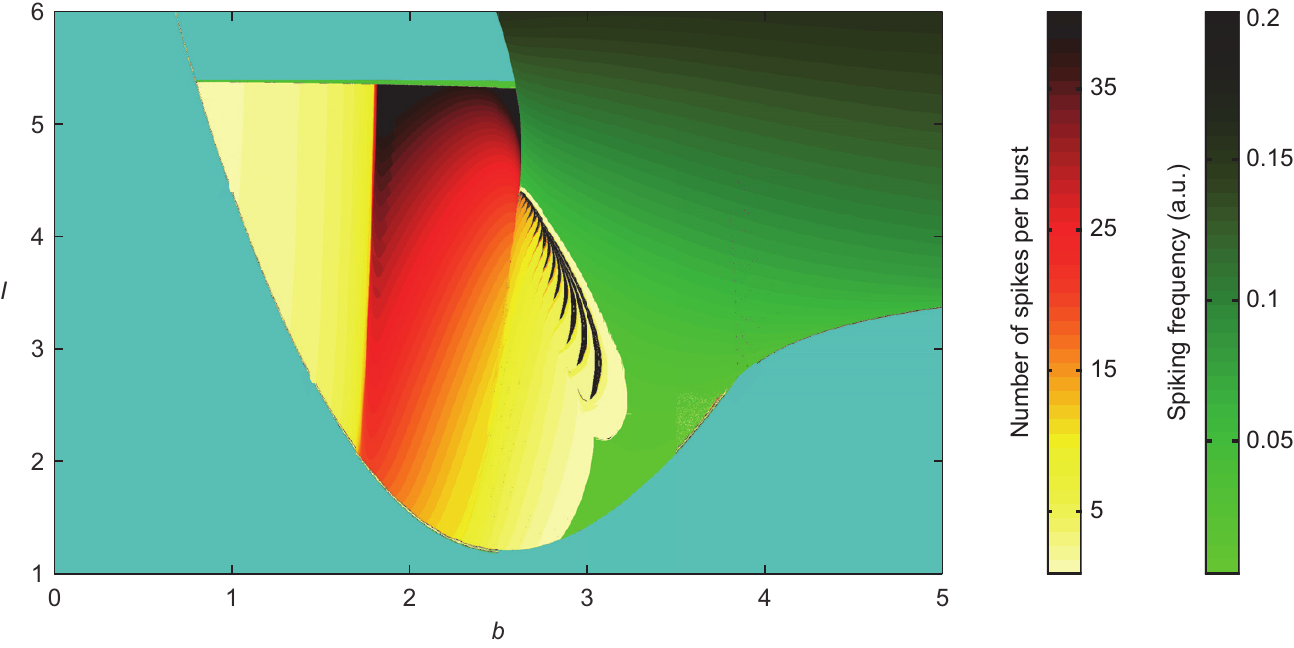}
\caption[Brute force bifurcation diagram of the HR model]
{Brute force bifurcation diagram of the HR model
\Ref{eq:HR}, for the set of parameter values~\Ref{eq:parvalues},
in the $(b,I)$-plane. See text for full explanation of the colour map.
}
\label{fig:HRbifwide}
\end{figure}

\InRefFig{fig:HRbifwide} shows a brute-force bifurcation diagram of the region of interest.
The bifurcation diagram indicates that the HR model is able to reproduce all
the dynamical behaviours indicated in the previous section. In particular, the colour code is as
follows:
\begin{itemize}
\item[\textbf{cyan}] represents quiescence;
\item[\textbf{green}] is for spiking, with darker tones corresponding to a higher steady-state
    firing frequency;
\item[\textbf{yellow to red colours}] are used to represent regular bursting (stable periodic
    orbits), more specifically yellow changes to red as the number of spikes per burst
    increases;
\item[\textbf{black}] represents irregular spiking, which, from a dynamical system's point of
    view, corresponds to chaotic behaviour.
\end{itemize}
See \cite{Storace:2008:Chaos} for an explanation of the techniques used. One of the weaknesses of
the method is that the presence of regions admitting coexisting asymptotic behaviours cannot be
directly inferred from the colour map.

\subsection{The regular-to-irregular bursting transition}

\begin{figure}[!t]
\centering\includegraphics[width=16cm]{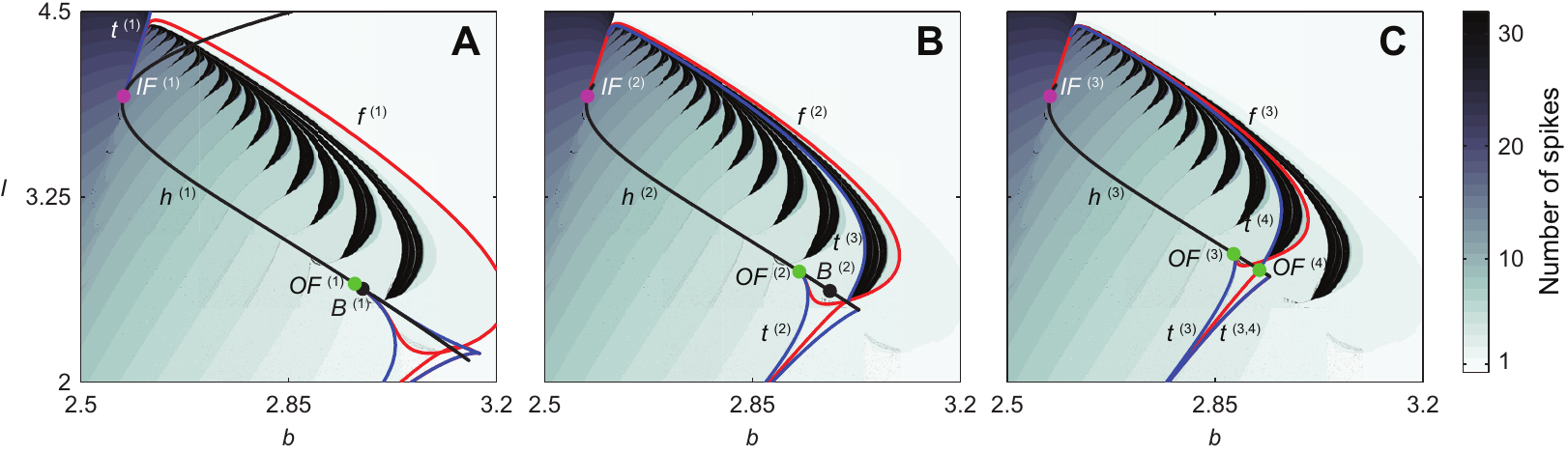}
\caption[Numerically computed bifurcation curves]{Numerically
  computed bifurcation curves showing
(in panels A, B and C respectively) bifurcations associated with
one, two and three bursts per period. See text for details.}
\label{fig:curves}
\end{figure}

In terms of bifurcation analysis, the area with the richest dynamics in \RefFig{fig:HRbifwide} occurs in
the region $b\in[2.5,3.2]$, $I\in[2,4.5]$.  Here we can observe lobe-shaped regions of irregular
bursting and stripe-shaped regions of regular bursting, with each successive stripe corresponding
to one extra spike per period. To try to understand the mechanism by which this transition from
regular to irregular bursting regions occurs, \RefFig{fig:curves} shows a zoom of the parameter region in
question with three different sets of bifurcation curves superimposed. These curves were computed
using numerical continuation in the software AUTO-07p\cite{Auto07PMan} and its extension HomCont
for the localisation of codimension-two homoclinic bifurcation points.

In the diagram, we have adopted the following colour and labelling code for the bifurcations curves
\begin{itemize}
\item folds of cycles (tangent bifurcations) are labelled $t$ and coloured blue;
\item period doublings (flips) are labelled $f$ and coloured red;
\item homoclinic bifurcations are labelled $h$ and coloured black.
\end{itemize}
Moreover, the superscript index indicates the approximate number of spikes per period of the limit
cycle (or homoclinic orbit) undergoing the bifurcation. So, for example, the label $f^{(1)}$
indicates a period doubling bifurcation curve involving a 1-spike cycle. Note that each of the
flips typically represents the first in an entire period-doubling cascade. Superscripts $(n,m)$
indicate that the cycle involved in the bifurcation undergoes a transition from $n$ to $m$ spikes
as is typical of bifurcations involved in the period adding mechanism. In addition we use letters
to distinguish distinct bifurcations of the same kind, e.g. $h^{(2)}$ and $h^{(2a)}$ will represent
different homoclinic orbits that have two spikes.

\begin{figure}[!t]
  \centerline{\includegraphics[width=12cm]{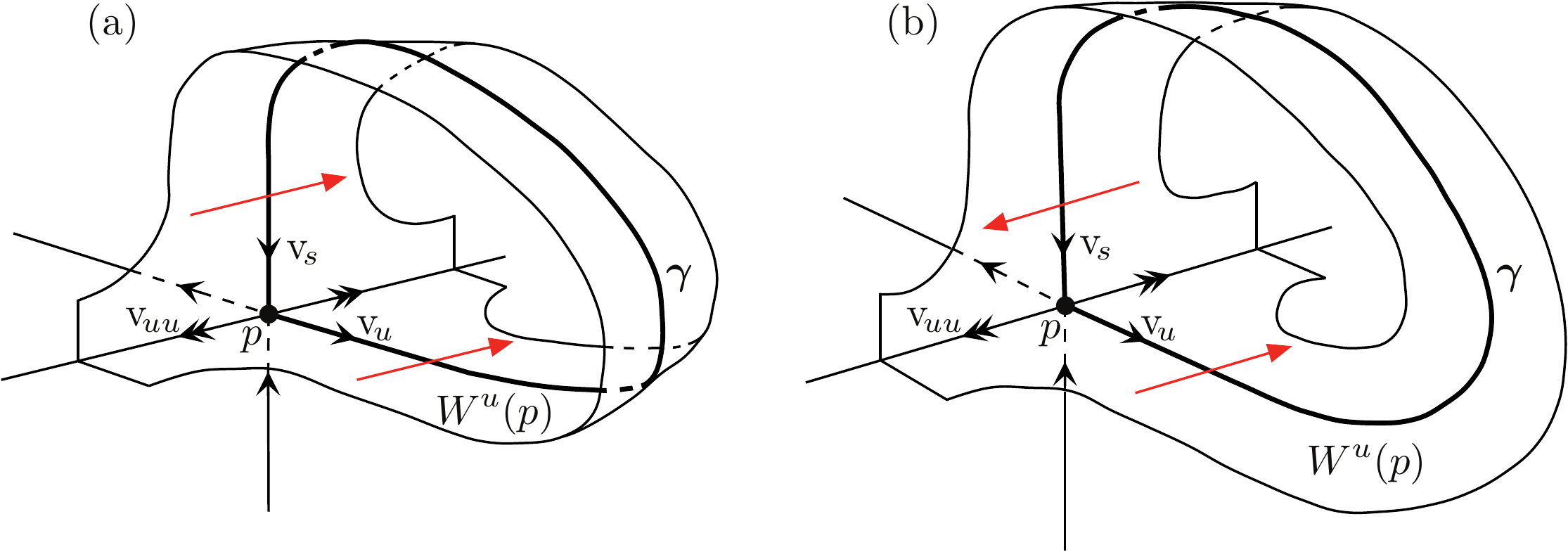}}
\caption{\label{fig:IF} Schematic representation of an inclination flip
bifurcation in the case (relevant to \RefEq{eq:HR})
of a saddle equilibrium in a three-dimensional vector field
with a two-dimensional unstable manifold.
The two panels show the same
homoclinic orbit before (panel (a)) and after (panel (b)) the bifurcation; note how the
unstable manifold changes from orientable to non-orientable
between panels (a) and (b), respectively.
The notation is as described in Sec.~\ref{sec:4} below.}
\end{figure}

In \RefFig{fig:HRbifwide} we have also identified several codimension-two homoclinic bifurcation
points. Specifically purple, green and black dots indicate respectively
  {\em inclination flip}, {\em orbit flip} and {\em Belyakov} points.
An inclination flip bifurcation represents a point along a curve of homoclinic orbits to a real
saddle at which the orientability of the global stable or unstable manifold changes, see
Fig.~\ref{fig:IF}. For information on the complex codimension-one curves that can emanate from the
codimension-two point see for example \cite{HomburgKrauskopf,Homburg:2010} and references therein. In particular, there are three
topologically distinct cases. An orbit flip occurs when the trajectory undergoing the homoclinic
orbit flips between the two components of the (weak) stable or unstable manifold. In the case of a
real saddle in three dimensions, the same three topological cases apply as for the inclination
flip, again see \cite{Sandstede,HomburgKrauskopf} and references therein. A Belyakov bifurcation
\cite{Belyakov:1974,Belyakov:1980,Belyakov:1984} occurs when the leading eigenvalues (closest to
the imaginary axis) of the saddle-point involved in the homoclinic orbit are double and undergo a
transition to a complex pair. The theory predicts the presence of several families (of infinite
cardinality) of bifurcation curves originating at these points and accumulating exponentially on
the homoclinic curve. See~\cite{Champneysetal:1994,Kuznetsov:2004,ShilnikovTuraevBook:2000} for
more details of the dynamics near codimension-two homoclinic bifurcations.
As we shall see in more detail shortly, these codimension-two points, and more besides, play a key
role in unfolding the regular to irregular bursting transitions.

Note that each of the three homoclinic bifurcation curves computed in \RefFig{fig:curves} actually
represents an approximate double cover of the same curve in parameter space, with the seeming
endpoint of the curve in fact representing a sharp $U$-turn. Therefore the fine structure of the
bifurcation curves is not apparent without looking at the particular shapes of the trajectories,
which will be elucidated in the following subsections. The structure of the homoclinic curve
$h^{(3)}$ and the associated local bifurcations of cycles depicted in panel C of
\RefFig{fig:curves} is similar to that relevant for all subsequent lobe-to-stripe transitions for
$k>3$.  Therefore the case $k=3$ will serve as an illustrative example in what follows. The cases
for $k=1$ and $k=2$ (depicted in panels $A$ and $B$ respectively) are special and will be dealt
with separately.

\begin{figure}[!t]
\centering
\includegraphics[width=10cm]{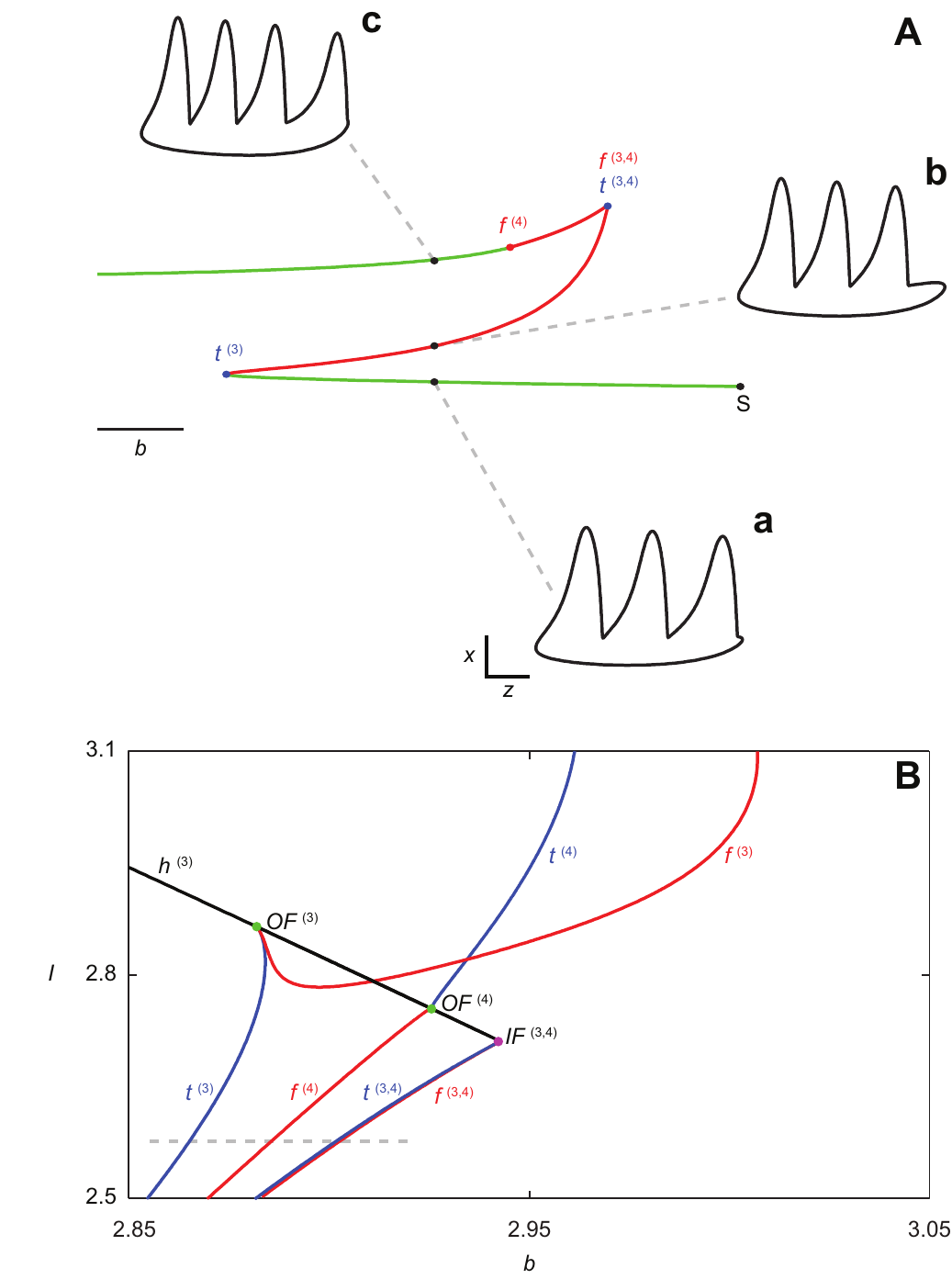}
\caption[Sketch of the period adding mechanism] {Sketch of the period
  adding mechanism and corresponding bifurcation curves.  In panel A
  the coloured traces indicate the maximum $z$ coordinate of the
  solution; colour encodes stability: green is stable, red is
  unstable. The trajectories a, b and c are projections on the $(z,x)$
  plane of the full three-dimensional solution.  Panel B shows the
  actual orbit flip points and the bifurcation curves that take part
  in the period adding mechanism. The continuation shown in panel A
  can be obtained by following, for example, the dashed grey line that
  crosses $t^{(3)}$, $f^{(4)}$, $t^{(3,4)}$ and $f^{(3,4)}$.  }
\label{fig:pasketch}
\end{figure}

Before proceeding with a more in-depth examination of the homoclinic bifurcations, it is worth
showing out how the local bifurcations of cycles that bifurcate {\em below} the homoclinic
bifurcation curves organise the boundaries of the stripe-shaped periodic bursting regions.
\InRefFig{fig:pasketch} shows in detail the bifurcations associated with the spike-adding boundary
between the 3-spike and 4-spike regular bursting regions. Note that the transition is
hysteretic; that is, there is a parameter window of bistability in which both 3- and 4-spike
regular bursting can be observed.

Panel A of the figure shows the nature of the transition under variation of $b$. Upon decreasing
the parameter from the point labelled $S$, the stable 3-spike cycle (labelled $a$) becomes unstable
through a fold $t^{(3)}$, leading to interval of unstable 3-spike cycles (such as the one labelled
$b$), until another fold $t^{(3,4)}$, after which it remains unstable. The branch then
re-stabilises at a period doubling bifurcation $f^{(4)}$ to form the four-spike cycle labelled $c$.
Upon further decrease of $b$, this stable 4-spike cycle will remain until a further bifurcation
$t^{(4)}$ (not shown in this sketch) and the whole process repeats for the 4- to 5-spike cycle
transition.  Thus from the point of view of a single limit cycle, the whole spike-adding process
from 1 to many can be thought of as a single {\em smooth} process.

Panel B of \RefFig{fig:pasketch} shows a zoom from \RefFig{fig:curves} C of the bifurcation curves that
are involved in the spike-adding mechanism. Note, in particular, that two of the crucial
codimension-one bifurcations involved $t^{(3)}$ and $f^{(4)}$ originate in the parameter plane from
two distinct orbit flips that lie on the apparent $h^{(3)}$ homoclinic bifurcation curve. The other
local bifurcation curves involved, $t^{(3,4)}$ and $f^{(3,4)}$, appear to originate from the
end-point of the homoclinic curve. As we shall show these bifurcation curves are actually caused by
an inclination flip that occurs at this apparent end point.

\subsection{The homoclinic curve $h^{(1)}$}

\begin{figure}[!t]
\centering
\includegraphics[width=10cm]{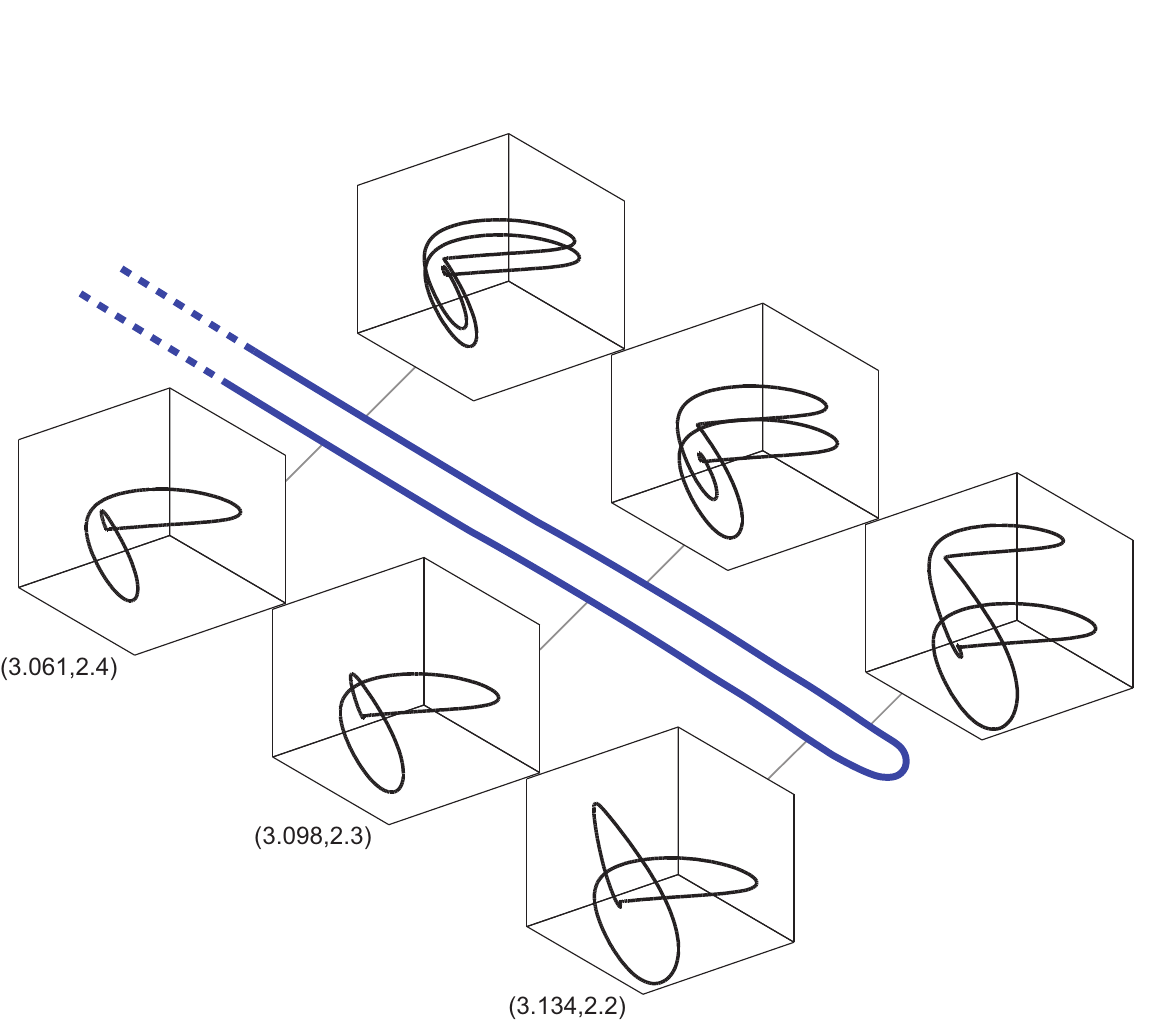}
\caption[Sketch of a U-shaped homoclinic bifurcation]
{Sketch of a U-shaped homoclinic bifurcation. The $(b,I)$
coordinates are reported
near the corresponding 1-spike trajectories on the lower branch of the homoclinic curve. The same values, with the chosen accuracy, hold for the 2-spike trajectories on the upper branch. Adapted from \cite{Linaro:2009},
reprinted with permission.}
\label{fig:Usketch}
\end{figure}

The first characteristic feature of each homoclinic curve $h^{(k)}$ is its U-shape, as
qualitatively sketched in \RefFig{fig:Usketch} for the lower part of the homoclinic curve $h^{(1)}$
(which curve is depicted only qualitatively). In fact, the U-turn is so sharp that it can be
detected only on a very small scale in the parameter space and on any wider scale, as in
\RefFig{fig:curves}A, the two branches appear almost as a double cover of the same curve. Note
that, as the U is traced, there is a transition between a homoclinic  with 1 spike to one with
2 spikes. Such sharp U-turns of homoclinic orbits have been observed in other systems, see
e.g.~\cite{Kuznetsovetal:2001,DeFeoetal:2000}, and are typically characterised by orbits gaining an
extra spike or pulse. As we shall see, similar U-turns on higher homoclinic branches $h^{(k)}$
cause transitions from $k$- to $(k+1)$-spike homoclinic cycles although, as we shall see, there are
important extra details for $k>3$, and the case $k=2$ is special.

\begin{figure}[!t]
\centering
\includegraphics[width=10cm]{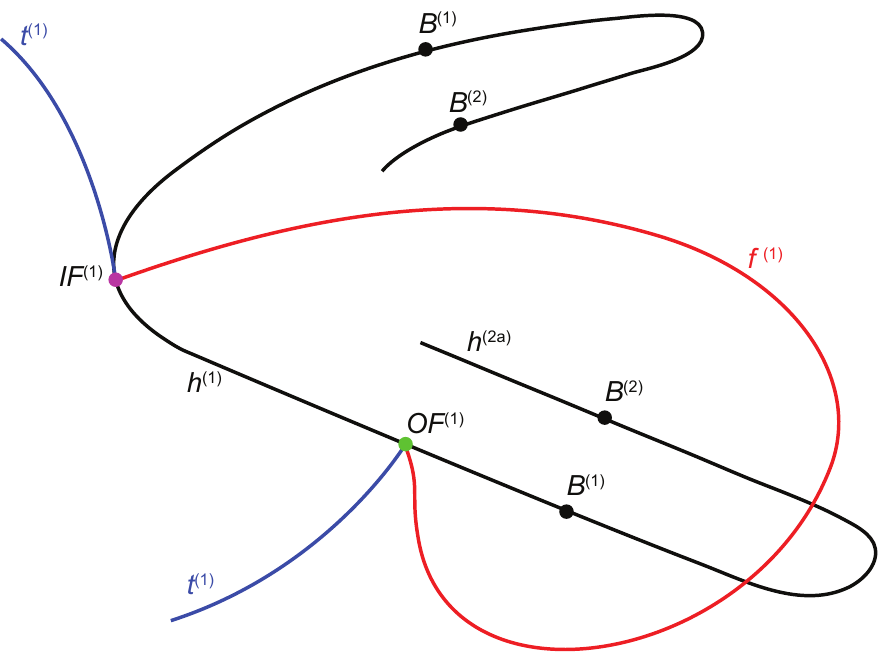}
\caption[Sketch of $h^{(1)}$]
{Sketch of $h^{(1)}$: note the presence of just one flip bifurcation ($f^{(1)}$) that connects
the inclination flip point $IF^{(1)}$ to the orbit flip point $OF^{(1)}$.}
\label{fig:sketch1}
\end{figure}

The homoclinic bifurcation curve $h^{(1)}$ is further sketched in \RefFig{fig:sketch1}. In this and
subsequent similar figures bifurcation curves are depicted in an exaggerated way in a pseudo
parameter plane in order to elucidate their topological features. The key feature here is an
inclination flip point labelled $IF^{(1)}$ that separates two portions of the homoclinic branch,
both U-shaped. We shall focus on the lower portion. The inclination flip in this case is of type B
according to the classification reported in \cite{HomburgKrauskopf} and \cite[Fig.7]{Oldeman2000};
there are single curves of period-doubling and fold of cycle bifurcations emanating from the
codimension-two point. Following the branch $h^{(1)}$ away from the inclination flip, we find an
orbit flip. Again, two other curves of local bifurcations emanate, a period doubling and a fold.
Note how the period-doubling bifurcation $f^{(1)}$ connects the two codimension-two points
$IF^{(1)}$ and $OF^{(1)}$, whereas the two fold of cycle curves labelled $t^{(1)}$ are distinct. We
also note the presence of Belyakov points, labelled as $B^{(1)}$ on the 1-spike branch $h^{(1)}$
and as $B^{(2)}$ on the 2-spike branch $h^{(2a)}$. Between these two points along the homoclinic
curve, the saddle equilibrium is actually a saddle focus (with complex eigenvalues). The real
curves are superimposed to the brute-force bifurcation diagram in panel A of \RefFig{fig:curves}.
With respect to the sketch, one further period-doubling and two fold of cycles (meeting at a cusp
point) curves are displayed (right-bottom corner). One of these fold of cycles is rooted at the
Belyakov point $B^{(1)}$.

The numerical continuation suffers convergence problems along the branches of 2-spike homoclinic
orbits as they return towards $IF^{(1)}$ and are depicted to end in ``mid air''. The eventual fate
of the multi-spike branches remains an open issue which we shall not address here, partly because
they do not seem to play any further role in the regular-to-irregular bursting transition of
interest in this paper.

\subsection{The homoclinic bifurcation curves $h^{(k)}$ and their degeneracies}
\begin{figure}[!t]
\centering
\includegraphics[width=10cm]{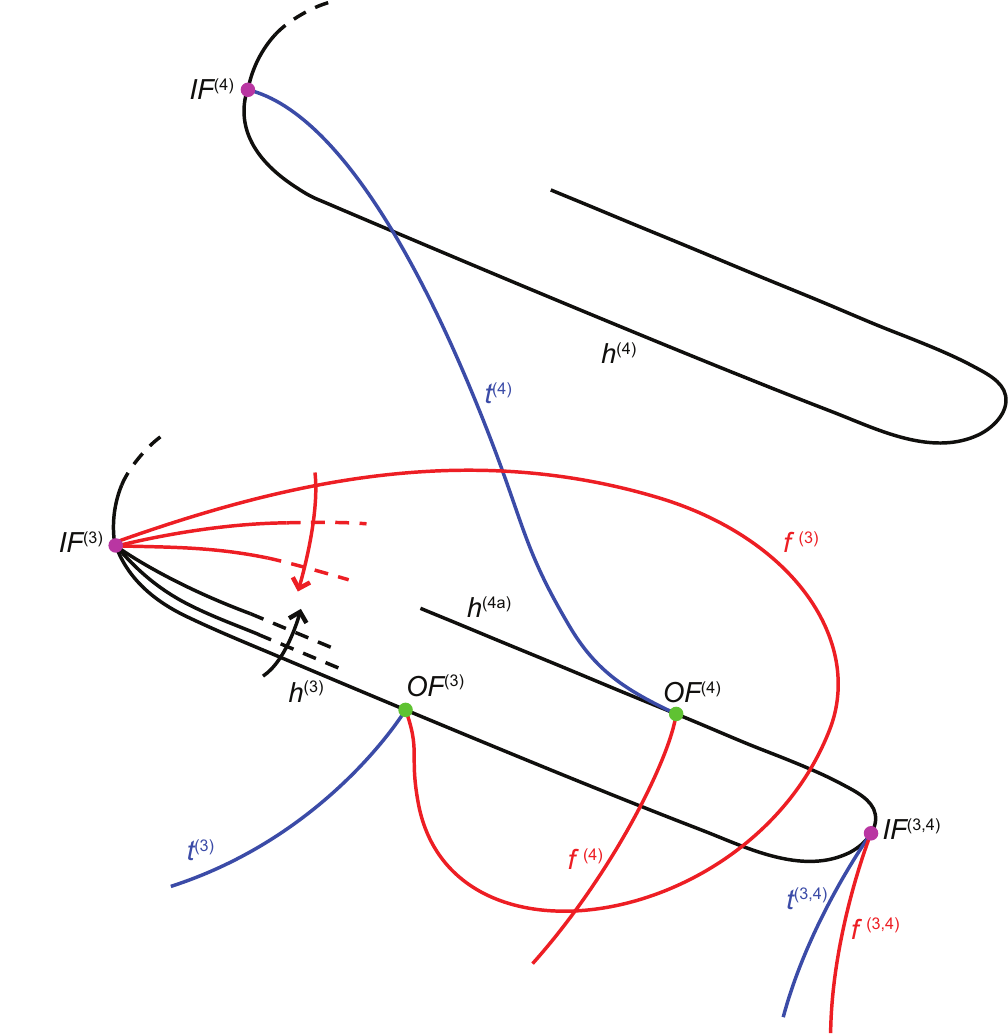}
\caption[Sketch of $h^{(k)}$ for $k\ge3$]
{Sketch of the bifurcation structure around the $k$-th homoclinic bifurcation $h^{(k)}$.
In this particular sketch, $k=3$, but the general sketch is valid for any $k\geq3$. For completeness, the sketch also shows the homoclinic curve for $k=4$.
This undergoes
the same sequence of bifurcations as $h^{(3)}$, which are not depicted. }
\label{fig:sketch3}
\end{figure}


The qualitative feature of the homoclinic bifurcation curve $h^{(3)}$, which is sketched in
\RefFig{fig:sketch3}, is valid for any $k \geq 3$. The fundamental difference with respect to
$h^{(1)}$ is the absence of any Belyakov point. This is due to the fact that the whole homoclinic
curve lies in a region of the parameters plane where the eigenvalues of the equilibrium involved in
the homoclinic trajectory are real. The curve $h^{(3)}$ emanates from an inclination flip point
$IF^{(3)}$, which appears to be distinct from $IF^{(1)}$ but occurs nearby in parameter space. One
distinction with the previous case is that the inclination flip is now of type C, which means that
an entire period-doubling cascade emanates, as do multiple-pulse homoclinic orbits for all periods
($h^{(3)}$,$h^{(6)}$,$\ldots$). These cascades are illustrated schematically via the sequence of
lines superimposed with curved arrows in \RefFig{fig:sketch3}. We do not focus here on these
additional bifurcations. Also, there are again computational difficulties with determining
precisely what happens to the branch on the ``far'' side of $IF^{(3)}$, or of the homoclinic curve
$h^{(4a)}$ as it returns towards the vicinity of $IF^{(3)}$. Instead we focus on the transitions
that occur close to U-turn as $h^{(3)}$ transitions into $h^{(4a)}$.

Each branch of the $U$-turn undergoes two orbit-flip bifurcations:
\begin{itemize}
\item $OF^{(3)}$, where the first bifurcation of the period-doubling cascade ends and meets a
    fold of cycles (respectively $f^{(3)}$ and $t^{(3)}$ in \RefFig{fig:sketch3}).
\item $OF^{(4)}$, where $t^{(4)}$ and $f^{(4)}$ are rooted: the former is connected with
    $IF^{(4)}$ on the primary homoclinic bifurcation of the subsequent homoclinic doubling
    cascade, whereas the latter takes part in the period adding process, as described in
\RefFig{fig:pasketch}.
\end{itemize}
The real curves are superimposed to the brute-force bifurcation diagram in panel C of
\RefFig{fig:curves}.

Very close to the tip of the U-shaped homoclinic bifurcation, there is an additional inclination
flip point, labelled $IF^{(3,4)}$ in \RefFig{fig:sketch3}. From this point, the two curves
$t^{(3,4)}$ $f^{(3,4)}$ are born, which are the additional bifurcations that take part in the
spike-adding process depicted in Fig.~\ref{fig:pasketch}. However, we have not been able to
numerically detect $IF^{(3,4)}$, due it would seem to the very sharp turn of the homoclinic curve,
but we can infer its presence as we now explain. Furthermore, the geometric analysis in
\RefSec{sec:4} shall provide more careful justification for the presence of this bifurcation.

\begin{figure}[!b]
\hspace{-4em}
\centering\includegraphics[width=0.8\textwidth]{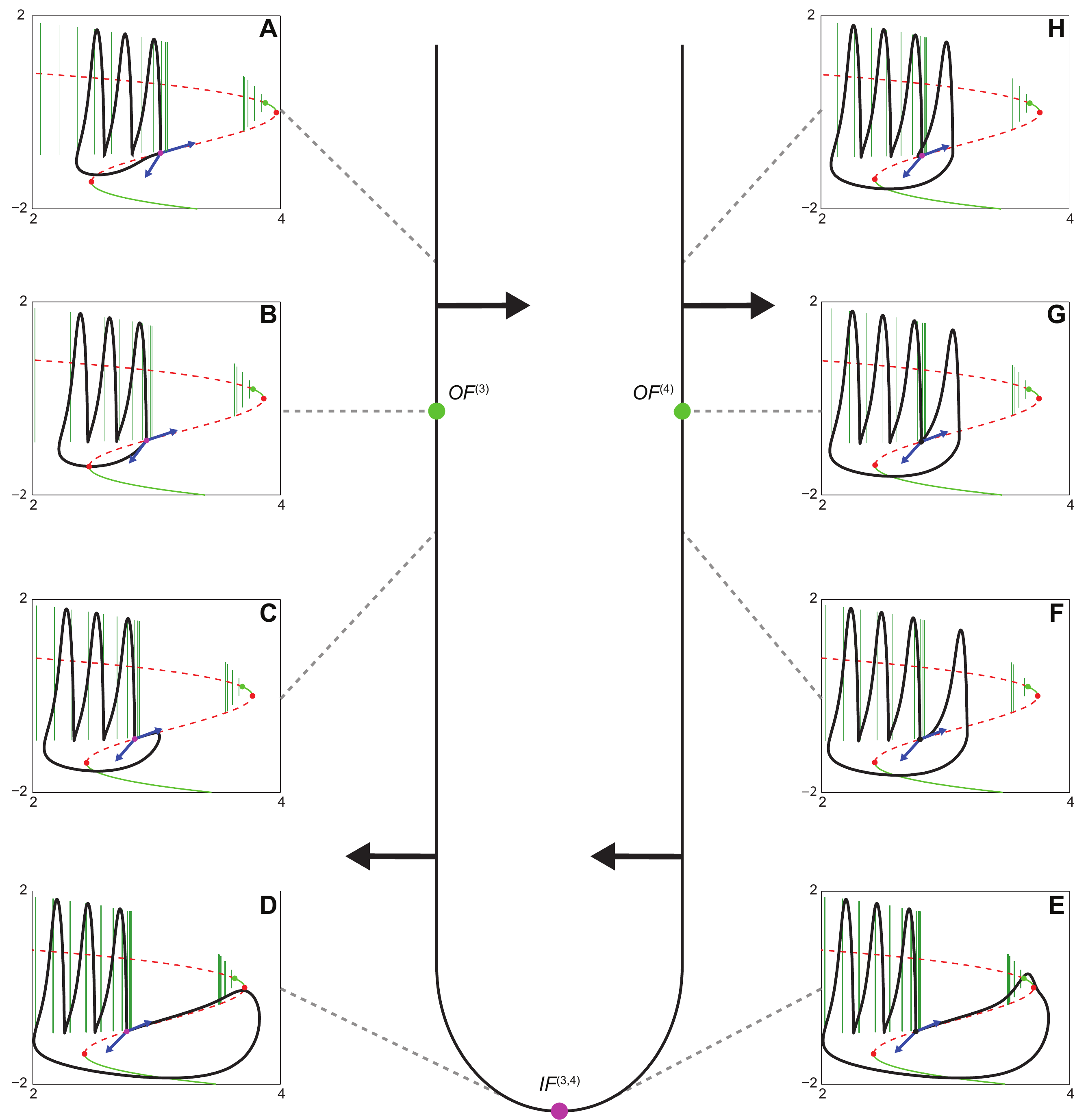}
\caption[Homoclinic trajectories along the bifurcation curve]
{
Representation of how the homoclinic trajectories change along the U-shaped homoclinic bifurcation
curve. Each panel contains the homoclinic orbit (thick black line) on the plane $(z,x)$ and the results of a bifurcation
analysis of the slow-fast subsystem of \RefEq{eq:HRfast} (thin coloured lines and dots); the blue
arrows are the unstable eigenvectors of the saddle node equilibrium. The arrows in the central panel, on the parameter plane,
indicate the direction of bifurcation of periodic orbits from the homoclinic bifurcation curve.
For a detailed description of each panel, see the main text.
}
\label{fig:IF_conjecture}
\end{figure}

\InRefFig{fig:IF_conjecture} shows more details of the orbits close to the $U$-turn, which provides
further evidence for the presence of the additional inclination flip $IF^{(3,4)}$. In this figure,
the central U-shaped curve represents the homoclinic bifurcation, and the eight surrounding panels
display the homoclinic trajectories (black thick lines) at significant points on the curve,
superimposed onto the bifurcation diagram of the fast subsystem (thin coloured lines
and points): these results are similar to those shown in \RefFig{fig:slowfast}, with the only
difference that the periodic solutions of the fast system do not constitute a unique ``funnel'',
but rather they are separated into two distinct sets, due to the presence of two homoclinic
bifurcations in the fast subsystem, at the coordinates where the periodic solutions accumulate.
Panels A and H are ``above'' $OF^{(3)}$ and $OF^{(4)}$, respectively, on opposite branches of the
homoclinic curve: in both panels, the homoclinic trajectory leaves the saddle node along the
leading unstable direction and returns along the only stable direction after 3 (panel A) or 4 turns
(panel H). Panels B and G correspond to the orbit flip points $OF^{(3)}$ and $OF^{(4)}$,
respectively: it can be clearly seen how the homoclinic trajectory leaves the saddle node along the
non-leading unstable direction. Again, the trajectory returns to the equilibrium point after 3
(panel B) or 4 spikes (panel G). Panels C-F are located between $OF^{(3)}$ and $OF^{(4)}$ and their
purpose is to illustrate the qualitative changes that the homoclinic trajectory undergoes between
the two orbit flip points and especially near the tip of the homoclinic curve, where we conjecture
the presence of the inclination flip point $IF^{(3,4)}$. In particular, in panels C and E it can be
observed how the homoclinic trajectory leaves the saddle node again along the leading unstable
direction, but this time in the opposite sense than in panels A and H. This makes the homoclinic
orbits sort of {\em canard cycles} that spend a large amount of time on the {\em unstable} part of
the slow manifold. In the $0<\mu\ll 1$-regime of the HR model it has been shown that canard
trajectories (not specifically homoclinic orbits) of this kind  exist for a wide range of
parameters, and they are known to be directly involved in the spike adding
mechanism~\cite{Desroches11,Guckenheimer,Terman:1991}. Finally, panels D and E are topologically
similar to panels C and F, with the only difference that, being so close to the tip of the
homoclinic curve, the canard orbits are \emph{maximal}: in particular, when the orbit goes past the
upper fold of equilibria in the fast subsystem, an additional turn is added to the trajectory,
which is the fundamental mechanism behind period adding in this and other models. Numerical
evidence shows that this happens exactly at the parameters values corresponding to the tip of the
homoclinic bifurcation curve.

The arrows in the central panel of \RefFig{fig:IF_conjecture} indicate the direction of bifurcation
of periodic orbits from the homoclinic bifurcation curve: the three points $OF^{(3)}$, $OF^{(4)}$
and $IF^{(3,4)}$ divide the homoclinic curve in four distinct regions. By going from one region to
the other, the direction of bifurcation of periodic orbits changes, due to the presence of the
orbit flip degeneracies and of the turning point at the tip of the U-shaped curve: this gives a
first, intuitive, indication that another degeneracy point where the orbit undergoes some switching
must be present. There are three such generic codimension-two points that lead to side-switching in
the case that the saddle point is a real saddle; orbit flip, inclination flip and resonant
eigenvalues. The latter occurs when $\mu_1 = -\lambda_1$, where $\mu_1$ is the stable eigenvalue of
the saddle point and $\lambda_1$ is the weakest unstable eigenvalue. We can easily check that the
eigenvalue condition is not satisfied, and we can rule out the presence of an orbit flip, since the
direction along which the trajectory leaves the saddle node does not change. Hence we are left only
with the possibility that the point at the tip is indeed an inclination flip. A very similar
structure has been found in \cite[Fig.19]{Champneys:2009} in another context.

\subsection{The special case $k=2$}
\begin{figure}[!t]
\centering
\includegraphics[width=10cm]{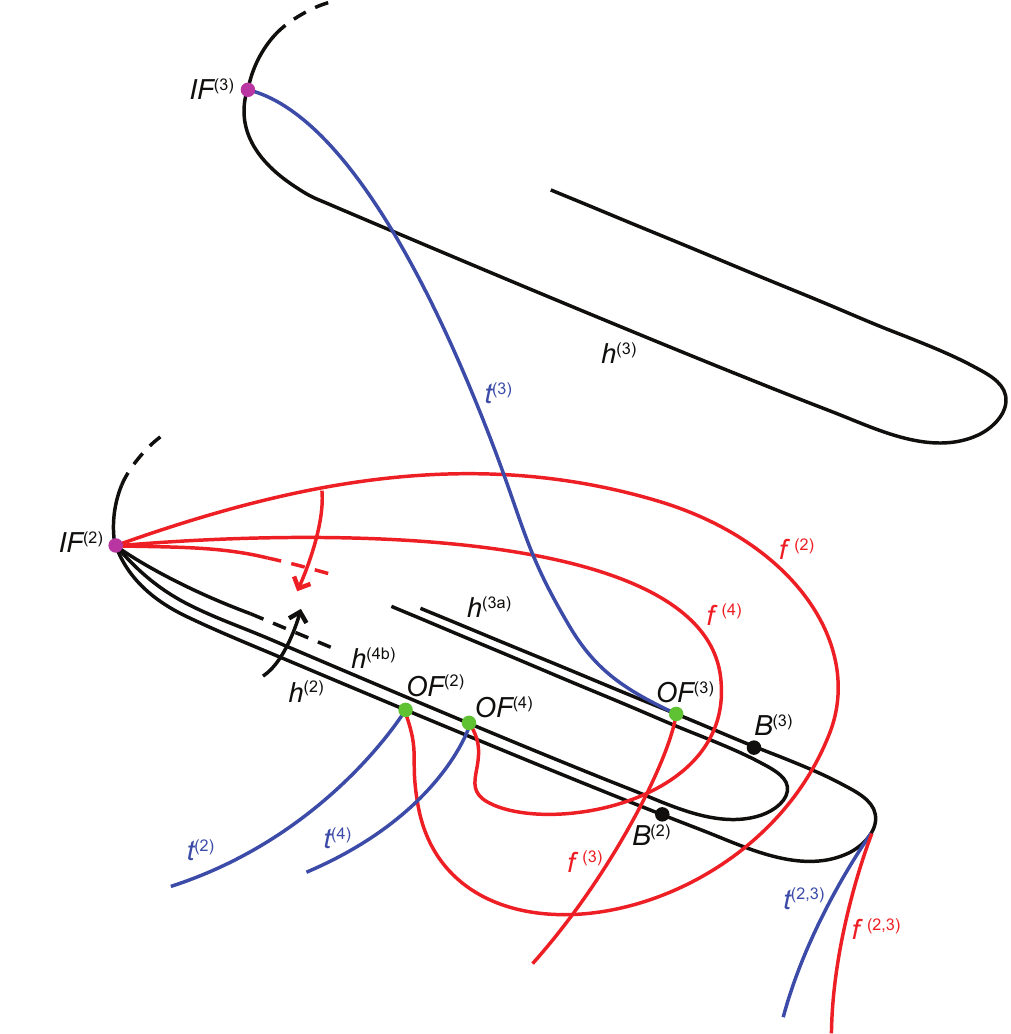}
\caption[Sketch of $h^{(2)}$]
{Sketch of $h^{(2)}$: the inclination flip point $IF^{(2)}$ gives birth to both a
homoclinic doubling and a period doubling cascade. Note the presence of two
Belyakov points, $B^{(1)}$ and $B^{(2)}$.}
\label{fig:sketch2}
\end{figure}

The homoclinic bifurcation curve $h^{(2)}$ is sketched in \RefFig{fig:sketch2} (the real curves
were superimposed on the brute-force bifurcation in panel B of \RefFig{fig:curves}). Here we also
compute a separate homoclinic curve $h^{(4b)}$ with four spikes that also comes out of the
inclination flip point $IF^{(2)}$ (which again appears to be distinct from $IF^{(1)}$ although
nearby to it in parameter space). This curve exists because the inclination flip is of type C and
is the first in an infinite sequence of the subsidiary homoclinic bifurcations that emanate from
the codimension-two point. Like in the general case for $k>2$, each homoclinic branch emanating
from the $IF$ has an orbit flip. The fold bifurcation $t^{(2)}$ is the one that is directly
involved in the spike-adding from 2 to 3 spikes similarly to what shown in \RefFig{fig:pasketch}
for the 3 to 4-spikes transition. A connection between the homoclinic curves $h^{(3a)}$ and
$h^{(3)}$ is provided by the fold of cycles $t^{(3)}$: this latter bifurcation terminates the
chaotic region that is born with the period doubling cascade that starts with $f^{(2)}$, as can be
seen in panel B of \RefFig{fig:curves}.

The curves $t^{(2,3)}$ and $f^{(2,3)}$ converge on the tip of the U-turn. Note that there can be no
inclination flip in this case, because between the two Belyakov points $B^{(2)}$ and $B^{(3)}$ the
equilibrium has complex eigenvalues.

%
\begin{figure}[!b]
  \centerline{\includegraphics[width=8cm]{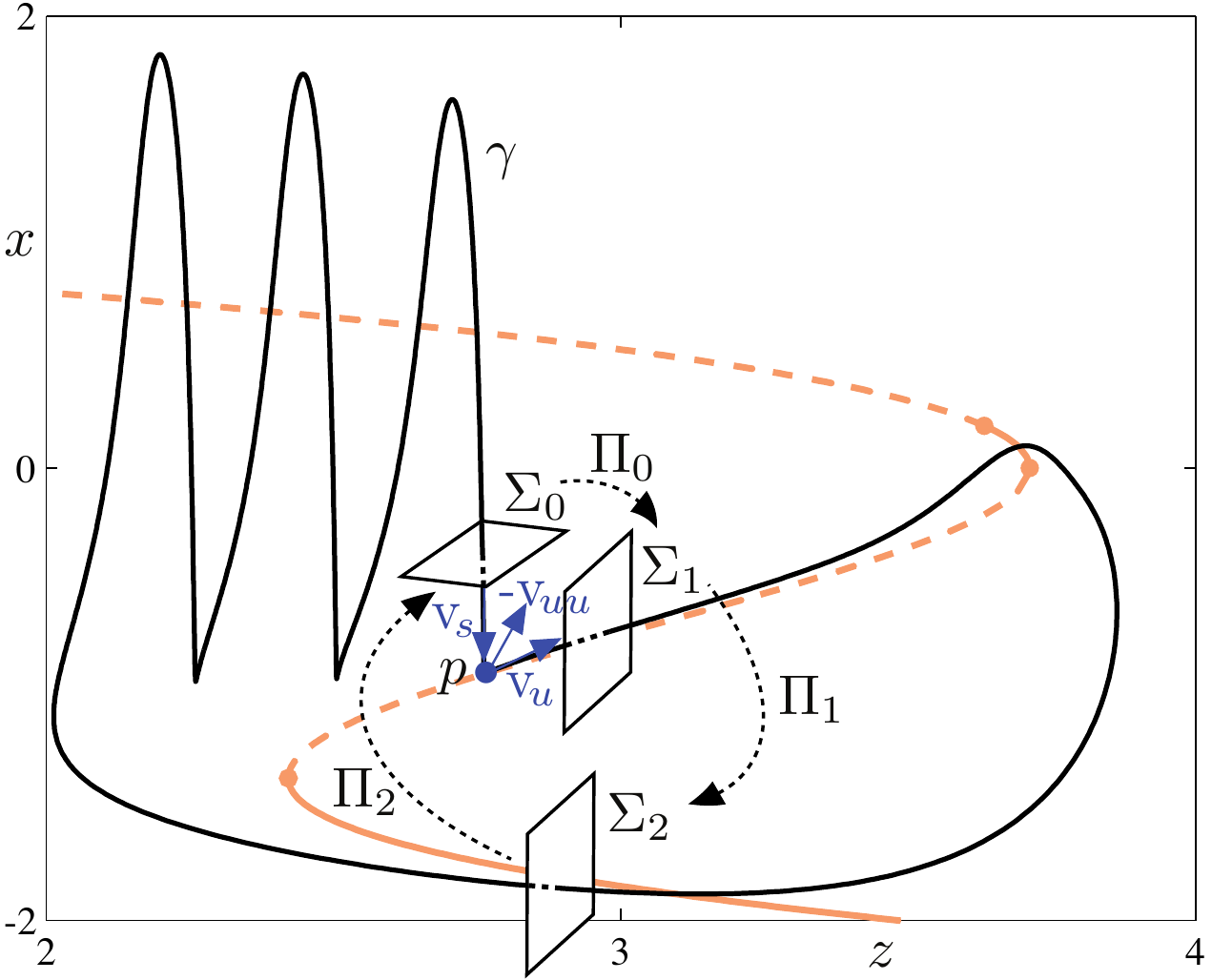}}
\caption{\label{fig:hom}
A 3-spike homoclinic orbit $\gamma$ of
system~(\ref{eq:HR}) projected onto the $(z,x)$-plane
  and superimposed onto the bifurcation diagram of the
  fast subsystem. The values of parameters $b$ and $I$ for this orbit
  correspond to the tip of the homoclinic curve computed
  in~\cite{Linaro11} and, hence, to the conjectured inclination flip
  bifurcation. We also show the saddle equilibrium $p$ together with
  its strong unstable, weak unstable and stable eigendirections
  $\mathrm{v}_{\!uu}$, $\mathrm{v}_{\!u}$ and $\mathrm{v}_{\!s}$,
  respectively. The three cross-sections $\Sigma_i$, $i=0, .., 2$,
  allow to construct a return map $\Pi$ from $\Sigma_1$ back to
  itself, in order to study the behaviour of nearby trajectories (for
  fixed $b$ and $I$ close to the transition of interest).}
\end{figure}
%
\section{Analysis of inclination flip due to fold in slow manifold}
\label{sec:4}

The purpose of this section is to show theoretically the presence of an inclination flip
codimension-two point at the sharp turning points of each of the curves $h^{(k)}$ with $k>2$.
Moreover, we aim to show that this process is a natural consequence of the sharp folding in the
curve of homoclinic orbits, and that this sharp turn is itself a consequence of the canard-related
transition of a $n$-spike homoclinic orbit into an $(n+1)$-spike homoclinic orbit. Furthermore, by
constructing an approximate return map around the critical homoclinic orbit, we are able to derive
asymptotic expressions for the curve of saddle-node of limit cycle bifurcations that emanates from
this codimension-two point.

The method of analysis is to construct the return map as a composition of approximate Poincar\'{e}
maps in a full neighbourhood of both parameter and phase space of the codimension-two point in
question; see Fig.~\ref{fig:hom}.  The analysis is general and can apply to any three-dimensional
system with the same generic features as the HR model. However, the key hypothesis has to be
justified numerically (in Subsection \ref{sec:if.1} which follows), namely that the forward image
of any smoothly parameterised set of trajectories that interacts transversely with the fold of the critical
manifold of the slow-fast system undergoes a sharp fold when viewed in any transverse Poincar\'{e}
section.  This assumption is formalised in the construction of the map $\Pi_2$ in Subsection
\ref{sec:if.2} below.

\subsection{The process of spike-adding}
\label{sec:if.1}

It is useful to examine in detail what happens to the trajectory of the homoclinic orbit as it
passes close to the sharp turning point in one of the loops of the loci of homoclinic orbits.
Consider \RefFig{fig:hom} which depicts just such an orbit that is undergoing a transition from
three to four spikes at parameter values $b = 2.9427488761, I = 2.7111448924$. We shall henceforth
refer to these as the critical parameter values.  Note that the nascent fourth spike forms via the
interaction of the unstable manifold with the fold point of the critical manifold (depicted by a
dashed red line). \InRefFig{fig:IF_conjecture} shows homoclinic orbits on the branch just before
(panel D) and just after (panel E) this critical codimension-two orbit.

\begin{figure}[!b]
\hspace{-3em}
\centering\includegraphics{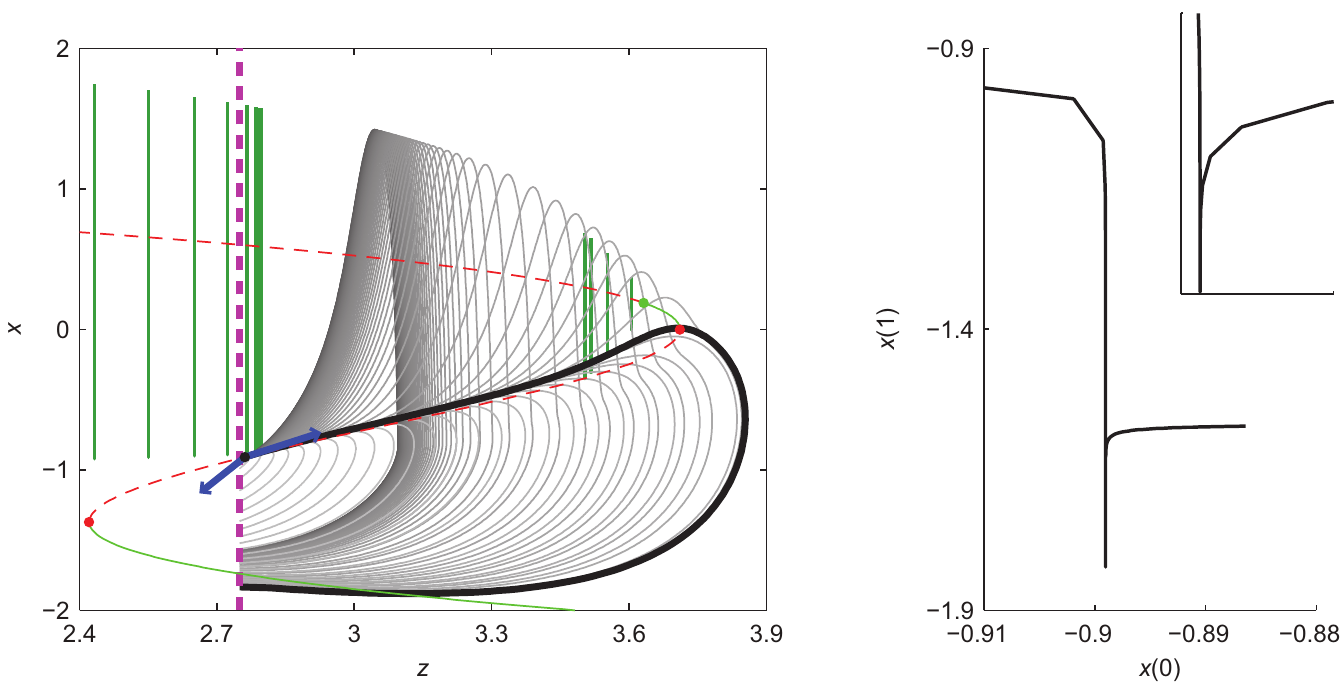}
\caption[Boundary Value Problem]{(Left panel) Trajectories in the unstable manifold at the
critical parameter values $b=2.9427488761$, $I=2.7111448924$. (Right panel) Approximate 1D  map
showing that the unstable manifold of the homoclinic trajectory is folded. For a detailed
description of each panel, see the main text.}
\label{fig:IF_BVP}
\end{figure}

\begin{figure}[!b]
  \centerline{\includegraphics[width=7.3cm]{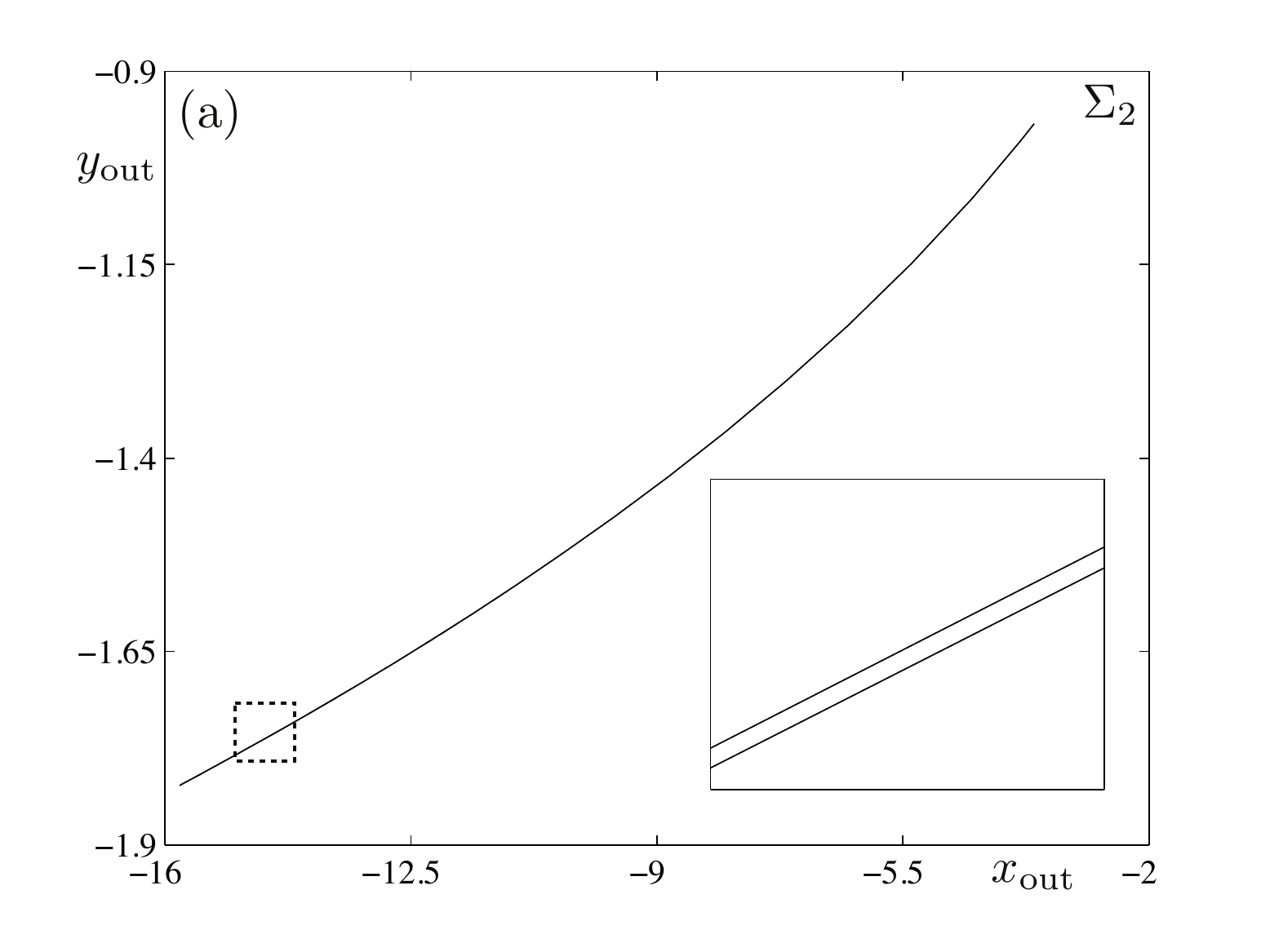} \quad
\includegraphics[width=6.5cm]{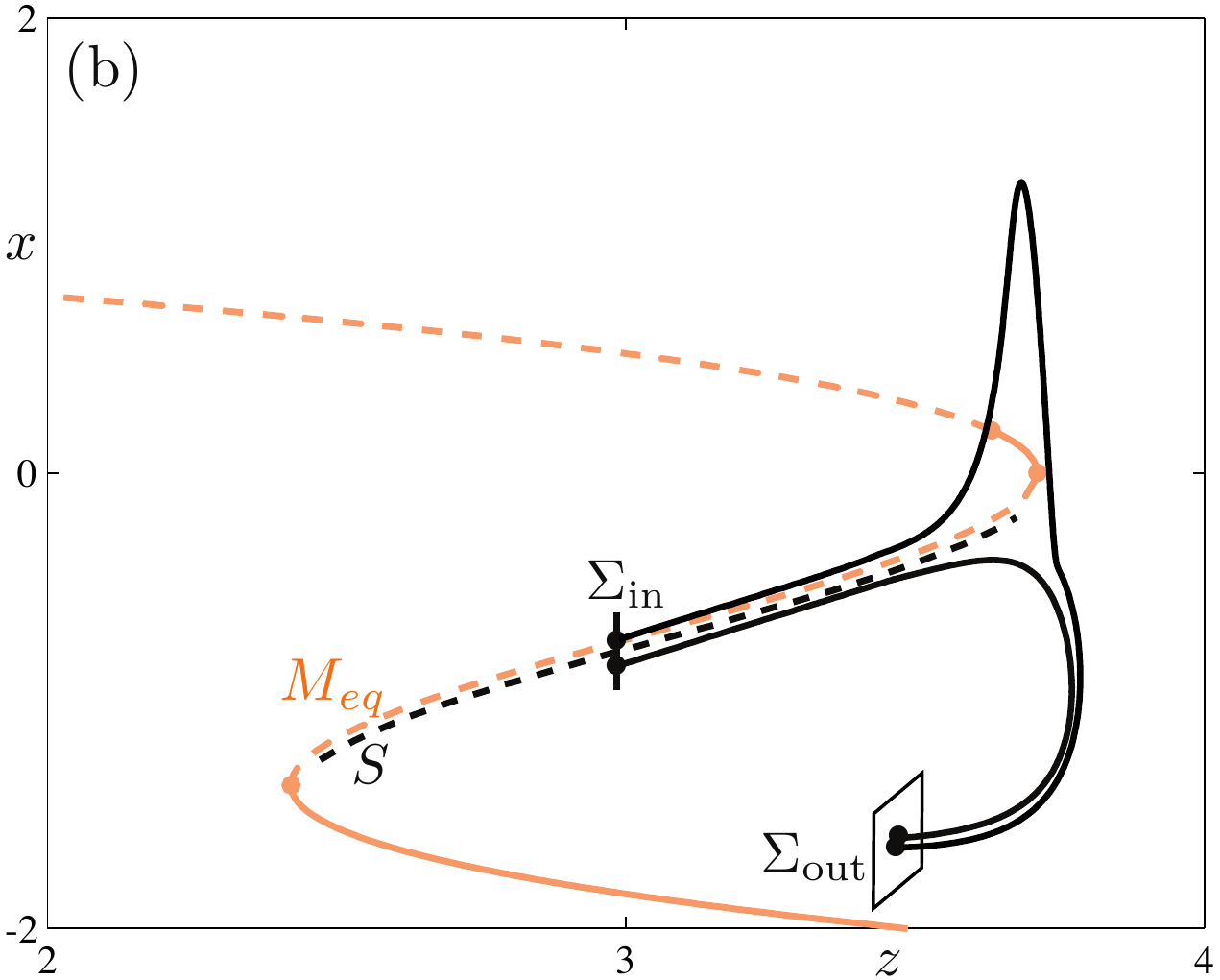}}
\caption{
(Panel (a)) Image of the initial conditions $U_1$ computed as in
\RefFig{fig:IF_BVP} projected onto the $x$ and $y$ co-ordinates of
the outgoing Poincar\'{e} section $z=2.75$.
(Panel (b))
Schematic representation of a general slow-fast system in 3D with a saddle-type
slow (Fenichel) manifold $S$ and an underlying critical manifold $M_{eq}$ that is folded.
Also shown are two orbit segments with very close initial conditions;
only one gets a twist when passing close to the upper fold of $M_{eq}$,
due to the relative position of its initial conditions with respect to $S$.}
\label{fig:numretmap}
\end{figure}

\InRefFig{fig:IF_BVP} and \RefFig{fig:numretmap} show the results of a numerical computation of a
portion of the unstable manifold of the saddle point at the critical parameter values. The map
shown in \RefFig{fig:IF_BVP} (right panel) was computed by variation of a transverse co-ordinate in
the unstable manifold close to the equilibrium point $p$ and computing until the first return to a
Poincar\'{e} section given by $z=2.75$. In particular the set $U_1$ of initial conditions chosen
was of the form
$$
U_1=\{(x,y,z)= p+ \eps \mathrm{v}_u + \theta \mathrm{v}_{uu} |
\mbox{ for } \theta \in
(-\eps, \eps) \}
$$
where $\eps=0.1$ was chosen to give a close approximation to the unstable manifold $W^u_{\rm
loc}(p)$ in a neighbourhood of the critical homoclinic orbit. Since the unstable manifold is an
invariant set, the theory predicts that trajectories that start on the manifold should remain on it
indefinitely: unfortunately, due to errors in the numerical integration of this slow-fast system
close to the critical manifold, such a result cannot be obtained with standard integration techniques.
However, it is possible to overcome this problem by resorting to \emph{continuation} techniques by
setting up a proper boundary value problem (BVP), where one of the parameters that are allowed to
vary is the integration time. This particular technique has been exploited, for example, in
\cite{Doedel:2006} to compute part of the manifold of the Lorenz system.

By solving the BVP with AUTO, we can obtain the results shown in \RefFig{fig:IF_BVP}. As in
previous figures, in the right panel the thin coloured lines represent
the bifurcation diagram of
the fast subsystem and the blue arrows are the unstable eigenvectors of the saddle node
equilibrium. The purple dashed line is the section that constitutes the terminating point of the
integrations. The thin grey lines are the integrations of the system obtained by varying the
parameter $\theta$ in the range $[-0.1,0.1]$; the thick black line is a piece of the homoclinic
trajectory that satisfies the boundary conditions and is used to start the continuation procedure
(which corresponds to $\theta\approx-0.001$). The left panel shows an approximate 1D  map of the
initial versus the final $x$-coordinate. It can be clearly seen that such an approximate map is
not invertible (see also the inset, which contains a zoom of the central part), \ie two distinct
initial conditions lead to the same final condition. This constitutes a further justification of
our conjecture, since it shows that the unstable manifold of the homoclinic trajectory is folded.

To show this folded manifold in more detail, we depict in
\RefFig{fig:numretmap} (a) the image of $U_1$ in the Poincar\'{e} section $z=2.75$. Note the folded
shape of the image of $U_1$.  We conjecture that this fold is a direct consequence of a portion of
the unstable manifold passing close to the fold point of the critical manifold $M_{eq}$. This
conjecture is confirmed by noting from the computation of the trajectories in question in
\RefFig{fig:IF_BVP} that the region of the sharp turning point in the image of $U_1$ corresponds to
the trajectories that pass the closest to the fold in $M_{eq}$ (actually since \RefFig{fig:IF_BVP} was
computed at the critical parameter values, the trajectory that corresponds to the closest point to the
fold is on the homoclinic orbit). Similar passage near such a fold of the critical manifold has
previously been found in the HR model and has previously been shown to underlie spike-adding at the
level of periodic orbits. It was first reported in~\cite{Terman:1991,Terman:1992}, where the author
focused on chaotic dynamics in between $n$-spike and ($n+1$)-spike orbits as well as the
disappearance of bursting upon parameter variation. This mechanism has been studied more recently
using the framework of slow-fast dynamical systems in~\cite{Guckenheimer}.

We show in~\RefFig{fig:numretmap} (b) a schematic representation of the dynamical behaviour
suggested by our numerical results. We depict a three-dimensional slow-fast system with all the
generic features of the HR model; two fast variables and, hence, a one-dimensional critical
manifold $M_{eq}$. The figure shows the projection onto the $(z,x)$-plane, where $z$ is slow and
$x$ is fast and we superimpose two segments of trajectories with initial conditions chosen to be
very close to one another and to $M_{eq}$. When $M_{eq}$ is cubic-shaped (i.e. with two fold
points) and when its middle branch is composed by saddle equilibria of the fast subsystem, then
away from the fold points this middle branch perturbs smoothly with respect to the small parameter
$\eps$ to a saddle slow (Fenichel) manifold $S$~\cite{Fenichel:1979}. In this configuration, one
can observe at the level of both transient and long-term dynamics nearby trajectories and
attractors that diverge from one another when one gains an extra twist as it passes close to the
upper fold of $M_{eq}$ whereas the other does not; see two such orbit segments in
\RefFig{fig:numretmap}.  This particular dynamical behaviour can be understood by further looking
at the underlying slow-fast structure of the problem. Indeed, the families of (un)stable manifolds
$W^{u,s}(p)$ of the saddle equilibria $p$ associated with the fast dynamics perturb smoothly to
stable and unstable manifolds $W^{u,s}(S)$ of the Fenichel manifold $S$~\cite{Fenichel:1979}. Then,
if two sets of initial conditions are taken close to the slow manifold $S$, on opposite sides of
the unstable manifold of $S$ will follow $S$ for some time until one jumps down and the other jumps
up, hence gaining an extra twist.

Thus, these numerical results provide strong justification that the process of spike adding is
caused by the portion of the trajectory of the homoclinic orbit that is closest in time to the
local unstable manifold passing close to the fold point of the slow manifold. In turn such a
passage causes a sharp fold in the forward image of the local unstable manifold. The aim of the
rest of this section is then to argue that this process causes a sharp turning point in parameter
space of the locus of homoclinic orbits and that there is necessarily an inclination-flip
bifurcation point there. Moreover, a fold curve of periodic orbits and a period doubling
bifurcation curve emanate from the inclination flip.

\subsection{Construction of Poincar\'{e} return map}
\label{sec:if.2}

Consider a sufficiently smooth three-dimensional vector field
$$
\dot{x} = f(x,\mu), \qquad x \in \R^3, \quad \mu \in \R^2,
$$
that has a saddle point $p$ with real eigenvalues $\lambda_{\!uu} > \lambda_{\!u} > 0
>\lambda_{\!s}$, with corresponding eigenvectors $\mathrm{v}_{\!uu}$, $\mathrm{v}_{\!u}$ and
$\mathrm{v}_{\!s}$. We assume for simplicity (after a parameter dependent change of co-ordinates if
necessary) that the location of and linearisation at $p$ is parameter independent. Suppose that, at a
critical codimension-two point $\mu=0$, a homoclinic orbit $\gamma(t)$ to $p$ exists that satisfies
certain non-degeneracy hypotheses:
\begin{description}
\item[(H1)] $\gamma(t) \to p$ as $t \to \pm \infty$
\item[(H2)] $\gamma(t)$ is tangent to $\mathrm{v}_{\!u}$ as $t \to -\infty$ and specifically
    approaches $p$ along the positive $\mathrm{v}_{\!u}$ direction.
\end{description}
We also suppose that the sign of $\mathrm{v}_{\!s}$ has been chosen so that $\gamma(t)$  approaches
$p$ along the positive $\mathrm{v}_{\!s}$ direction as $t \to +\infty$.
\begin{description}
\item[(H3)] The map $\Pi_2$  (defined below) is degenerate such that it can be described by the
    given quadratic form to leading order.
\item[(H4)] The parameter $\mu$ unfolds this codimension-two singularity in a generic way.
    (Specific choices for $\mu=(\mu_1,\mu_2)$ are defined below.)
\end{description}

\begin{figure}[!b]
  \centerline{\includegraphics[width=8cm]{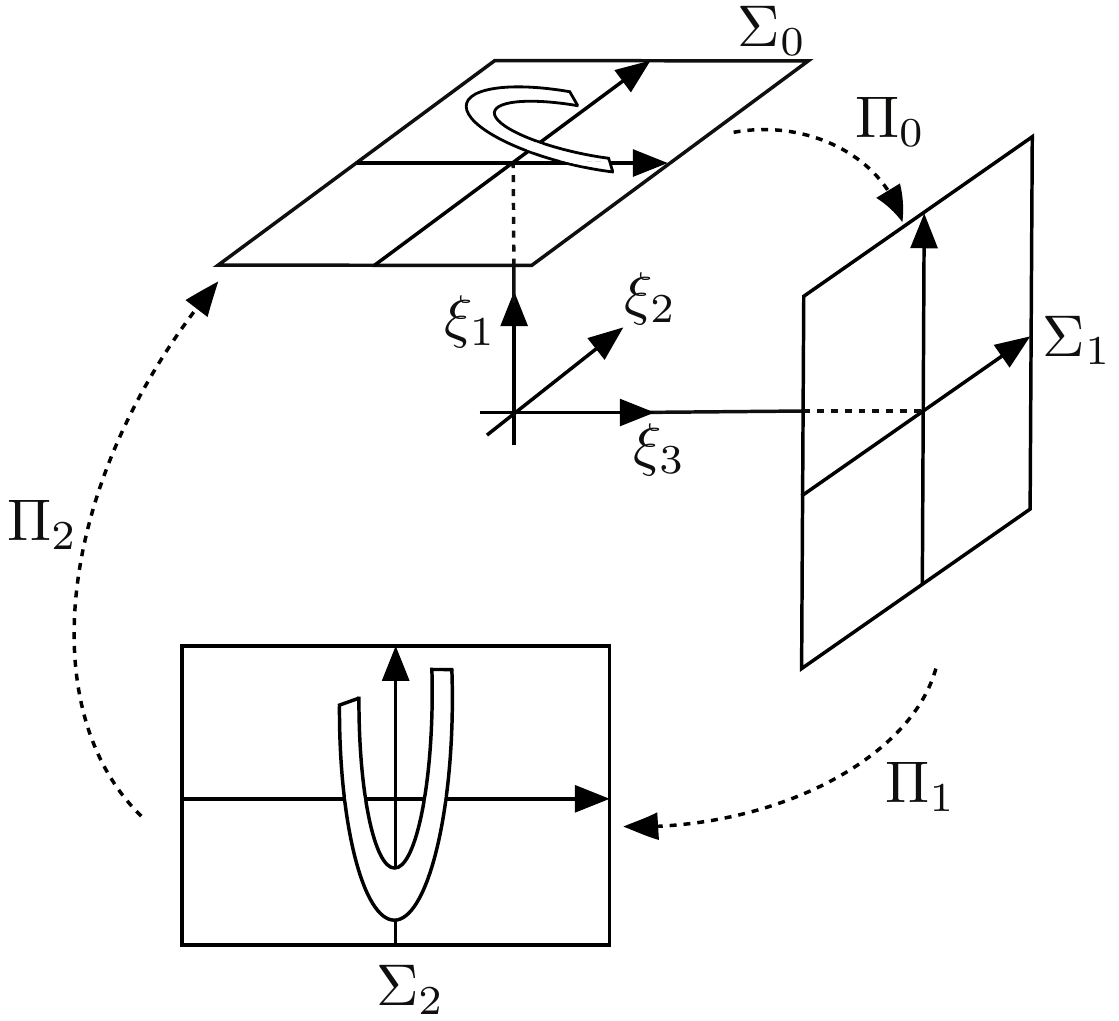}}
\caption{\label{fig:gensys}
Poincar{\'e} sections $\Sigma_0$, $\Sigma_1$, and $\Sigma_2$ for the study of the inclination flip
bifurcation in a general three-dimensional system.}

\end{figure}

We begin the analysis by considering three separate Poincar\'{e} sections $\Sigma_0$, $\Sigma_1$,
and $\Sigma_2$ as depicted in Fig.~\ref{fig:hom} for the HR model and in Fig.~\ref{fig:gensys}
for a general system.
The cross-sections $\Sigma_0$ and $\Sigma_1$ are defined in terms of local coordinates $(\xi_1,
\xi_2, \xi_3)$ corresponding to projection along the three-dimensional basis ($\mathrm{v}_{\!s}$,
$\mathrm{v}_{\!u}$, $\mathrm{v}_{\!uu}$). Specifically, let
$$
\Sigma_0=\{(\xi_1,\xi_2,\xi_3)| \xi_1=\eps\}, \qquad
\Sigma_1=\{(\xi_1,\xi_2,\xi_3)| \xi_2=\eps\}
$$
for $0<\eps \ll 1$.

\begin{figure}[!t]
  \centerline{\includegraphics[width=10cm]{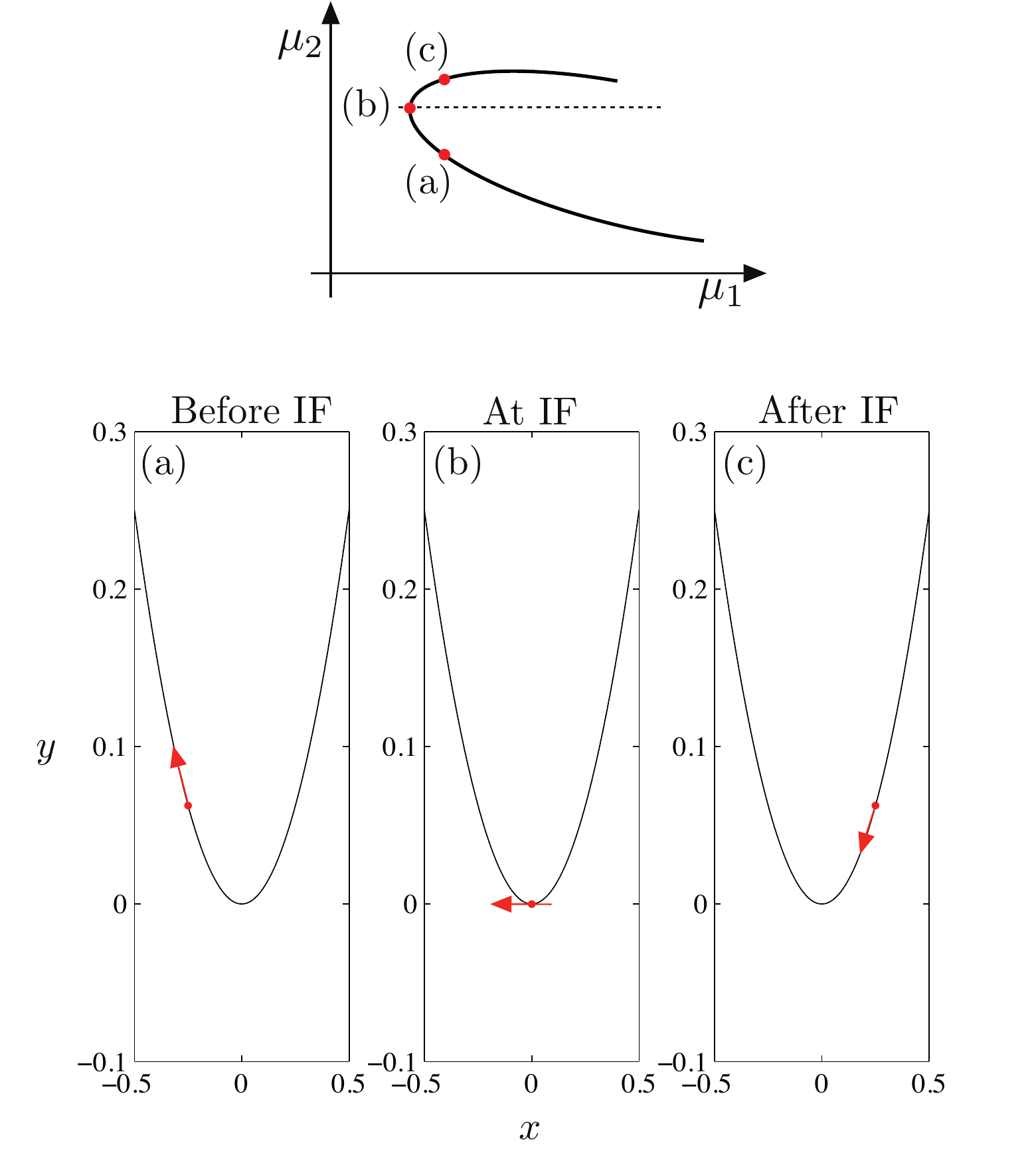}}
\caption{\label{fig:schretmap}Schematic representation of the return map computed for this system.}
\end{figure}

The section $\Sigma_2$ is chosen to be transverse to the flow at a point $\gamma(0)$ along the
critical homoclinic orbit, at an $O(1)$ distance from $p$. Let local co-ordinates $(\eta_1,\eta_2)$
be chosen within $\Sigma_2$ such that $\gamma(0)$ is at the origin and the tangent vector to
$W^u(p)\cup \Sigma_2$ at $\gamma_0$ lies along the $\eta_1$ axis. Furthermore, after a parameter
dependent change of co-ordinates if necessary, we shall suppose that the flow from $\Sigma_1$ to
$\Sigma_2$ is independent of the unfolding parameters $\mu$. A convenient choice of unfolding
parameters is to assume that the intersection between $\Sigma_2$ and the component of the
one-dimensional stable manifold $W^s(p)$ that corresponds to $\gamma(0)$ when $\mu=0$ is precisely
given by $(\eta_1,\eta_2)=(\mu_1, \mu_2)$. See Fig.~\ref{fig:schretmap}.

We are now in the position to define leading-order Poincar\'{e} maps obtained by following
trajectories between each of these Poincar\'{e} sections. The local map from $\Pi_0: \Sigma_0 \to
\Sigma_1$ can be obtained by solving the linear equations in a neighbourhood of $p$.  It is most
useful in what follows to instead deal with $\Pi_0^{-1}:\Sigma_1 \to \Sigma_0$. Specifically, to
leading order we obtain
$$
\Pi_0^{-1} : \begin{pmatrix} \xi_1 \\ \eps \\ \xi_3
\end{pmatrix}
\mapsto \begin{pmatrix} \eps \\ \xi_2 \\ \xi_3
\end{pmatrix} =
\begin{pmatrix}
\eps \\
K_1 \xi_1^{\Delta_1} \\
K_2 \xi_3 \xi_1^{\Delta_2}
\end{pmatrix},
$$
where
$$
0< \Delta_2 =
\frac{\minusone\lambda_{\!uu}}{\lambda_{\!s}}<\Delta_1 = \frac{\minusone\lambda_{\!u}}{\lambda_{\!s}} , \quad
K_1 = \eps^{1-\Delta_1} \quad \mbox{and } K_2 = \eps^{-\Delta_2}
$$
Hypothesis (H3) can now be encapsulated  in the leading-order expression for the Poincar\'{e} map
$\Pi_1 :\Sigma_1 \to \Sigma_2$. We construct this map in two stages. First consider the image of
$W^u(p)$
$$
\begin{pmatrix}
\eta_1 \\ \eta_2
\end{pmatrix}
=
\begin{pmatrix}
\xi_3 \\
\beta  \xi_3^2
\end{pmatrix}
$$
where the unit coefficient of the $\eta_1$-term is chosen without loss of generality. Also, the
$\eta_2$ co-ordinate is chosen so that $\beta>0$. Moreover, the assumption (H3) of a sharp fold in
the image of $W^u(p)$ implies
\begin{equation}
\beta \eps \gg 1
\label{eq:betascale}
\end{equation}
Thus the leading-order expression for the unit tangent vector to $W^u(p) \cap \Sigma$ is
\begin{equation}
\label{eq:D}
\tau(\xi_3) = \begin{pmatrix}
D(\xi_3) \\
2\beta \xi_3 D(\xi_3)
\end{pmatrix}
, \quad \mbox{where} \quad
\quad D(\xi)=\frac{1}{\sqrt{1+4 \beta^2 \xi^2}}
\end{equation}
from which we obtain that the unit normal (in the sense of positive $\xi_1$ co-ordinate) is
$$
\tau^\perp (\xi_3) =
\begin{pmatrix}
-2\beta \xi_3 D(\xi_3) \\
D(\xi_3)
\end{pmatrix}
$$
Hence, the leading-order expression for the full map $\Pi_1$ can be written as
$$
\Pi_1 : \begin{pmatrix}
\xi_1 \\ \eps \\ \xi_3
\end{pmatrix}
\mapsto
\begin{pmatrix}
\eta_1 \\
\eta_2
\end{pmatrix}
=
\begin{pmatrix}
\xi_3-2\beta \xi_1 \xi_3 D(\xi_3) \\
\beta \xi_3^2 + \xi_1 D(\xi_3)
\end{pmatrix}
$$

Finally, we suppose that the mapping $\Pi_2: \Sigma_2 \to \Sigma_0$ is a diffeomorphism that can be
expressed to leading-order by its linear terms.
$$
\Pi_2:
\begin{pmatrix}
\eta_1 \\
\eta_2
\end{pmatrix}
\mapsto
\begin{pmatrix}
\xi_2 \\
\xi_3 \\
\end{pmatrix}
=
B
\begin{pmatrix}
\eta_1 - \mu_1 \\
\eta_2 - \mu_2
\end{pmatrix},
$$
where $B = \{b_{ij}\}_{i,j=1,2}$ can be assumed generically to be an invertible matrix with all
elements nonzero.

\subsection{The inclination-flip bifurcation}

In the context of the example system in question, an inclination flip is a codimension-two
bifurcation that occurs when a path of homoclinic orbits to $p$ undergoes a change in orientation.

From the construction above, the locus of homoclinic orbits to $p$ in the $\mu$-plane is given to
leading order by \beq \mu_2 = \beta \mu_1^2 \label{eq:homaxis} \eeq which describes a sharp folded
curve pointing along the positive $\mu_2$-axis. For parameter values within this curve, the
twistedness of the unstable manifold along the homoclinic loop $\gamma$ can be computed by
following the tangent vector to the stable manifold around the homoclinic orbit $\gamma(t)$.

Let $(\mu_1,\mu_2)$ be a point within
the homoclinic locus given by \eq{eq:homaxis} and consider such a tangent vector with initial
condition in the positive $\mathrm{v}_{\!uu}$ direction within $\Sigma_1$.  By construction, the
image of this initial condition under $\Pi_1$ is the vector $\tau(\mu_1)$ defined above (Fig.~\ref{fig:hom}). The image
of $\tau(\mu_1)$ under $\Pi_2$ is then
$$
\hat{\tau}(\mu_1) :=
\begin{pmatrix}
(b_{11} + 2\beta b_{12} \mu_1 ) D(\mu_1)    \\
(b_{21} + 2\beta b_{22} \mu_1) D(\mu_1)
\end{pmatrix}
$$
Consider $\hat{\tau}$ for $\mu_1=\eps$. To leading-order we find
$$
\hat{\tau}(\eps) =
\begin{pmatrix}
b_{12}  \\
b_{22}
\end{pmatrix}
$$
in which we have used the form of $D$ defined above (\ref{eq:D}) and
the scaling (\ref{eq:betascale}). Now, under the non-degeneracy hypothesis that $b_{22} \neq 0$, as
$t \to \infty$, the tangent vector will tend to $\mbox{sign} (b_{22}) \mathrm{v}_{\!uu}$.

A similar argument shows that $\hat{\tau}$ for $\mu_1=-\eps$ along the homoclinic locus maps to
leading order to
$$
\hat{\tau}(\eps) =
\begin{pmatrix}
-b_{12} \\
-b_{22}
\end{pmatrix}
$$
in $\Sigma_0$. In turn, this vector tends to $-\mbox{sign} (b_{22}) \mathrm{v}_{\!uu}$ as $t \to
\infty$.

Hence we have shown that the tangent vector to the homoclinic orbit which is in the positive
$\mathrm{v}_{\!uu}$ component as $t \to -\infty$, flips its $\mathrm{v}_{\!uu}$ component as $t \to
+\infty$, for $\mu$ varying along the homoclinic locus \eq{eq:homaxis} between $\mu_1=\eps$ and
$\mu_1=-\eps$.  This shows that
there must be an (at least one) inclination flip somewhere in between. In other words, there must
be an orbit flip close to the sharp fold in the homoclinic locus.

\subsection{Unfolding the dynamics near the inclination flip}

Here we will extend the analysis to provide a local asymptotic prediction of the bifurcations of
periodic orbits that emanate from the inclination flip point. To this end we look for fixed points
of the return map
$$
\Pi_2 \circ \Pi_1 \circ \Pi_0  :\Sigma_0 \to \Sigma_0.
$$
In fact, it is most convenient to consider the Poincar\'{e} section $\Sigma_2$ and seek a condition
for a fixed point in the form
$$
\Pi_2^{-1} \circ \Pi_0^{-1} (\xi_1, \xi_3)^T = \Pi_1 (\xi_1, \xi_3)^T.
$$
To this end we find
$$
\begin{pmatrix} \mu_1 \\ \mu_2 \end{pmatrix}
+ \hat{B}
\begin{pmatrix} K_1 \xi_1^{\Delta_1} \\
                K_2 \xi_3 \xi_1^{\Delta_2}
\end{pmatrix} =
\begin{pmatrix}
\xi_3(1- 2 \beta \xi_1 D(\xi_3) ) \\
\beta \xi_3^2 + \xi_1 D(\xi_3),
\end{pmatrix}
$$
where
$$
\hat{B} = B^{-1} = \frac{1}{\mbox{ det}(B)}
\begin{pmatrix}
b_{22} & -b_{12} \\
-b_{21} & b_{11}
\end{pmatrix}
$$

We now need to analyse these fixed point equations and find fold and flip bifurcations. We suppose
we can do a rescaling so that $B=\textrm{Id}$. Then the equations for the fixed points of the
return map read
\begin{eqnarray*}
\mu_1+K_1\xi_1^{\Delta_1}          &=& \xi_3-2\beta\xi_1\xi_3 D(\xi_3)\\
\mu_2+K_2\xi_3\xi_1^{\Delta_2} &=& \beta\xi_3^2+\xi_1 D(\xi_3).
\end{eqnarray*}
Given that $\beta$ is assumed to be large ($\beta\eps\gg 1$), we make the following approximation
for $\xi_3$ such that it is at least of order 1, that is, ``non-small".
\beq D(\xi_3):=\frac{1}{\sqrt{1+4\beta^2\xi_3^2}}\approx\frac{1}{2\beta|\xi_3|}. \eeq
The fixed point equations (multiplying the second one by $\xi_3$) then reduce to
\begin{eqnarray*}
\mu_1+K_1\xi_1^{\Delta_1}          &=& \xi_3\mp\xi_1\\
\mu_2\xi_3+K_2\xi_1^{\Delta_2}\xi_3^2 &=& \beta\xi_3^3\pm\frac{\xi_1}{2\beta}.
\end{eqnarray*}

\begin{figure}[!t]
  \centerline{\includegraphics[width=10cm]{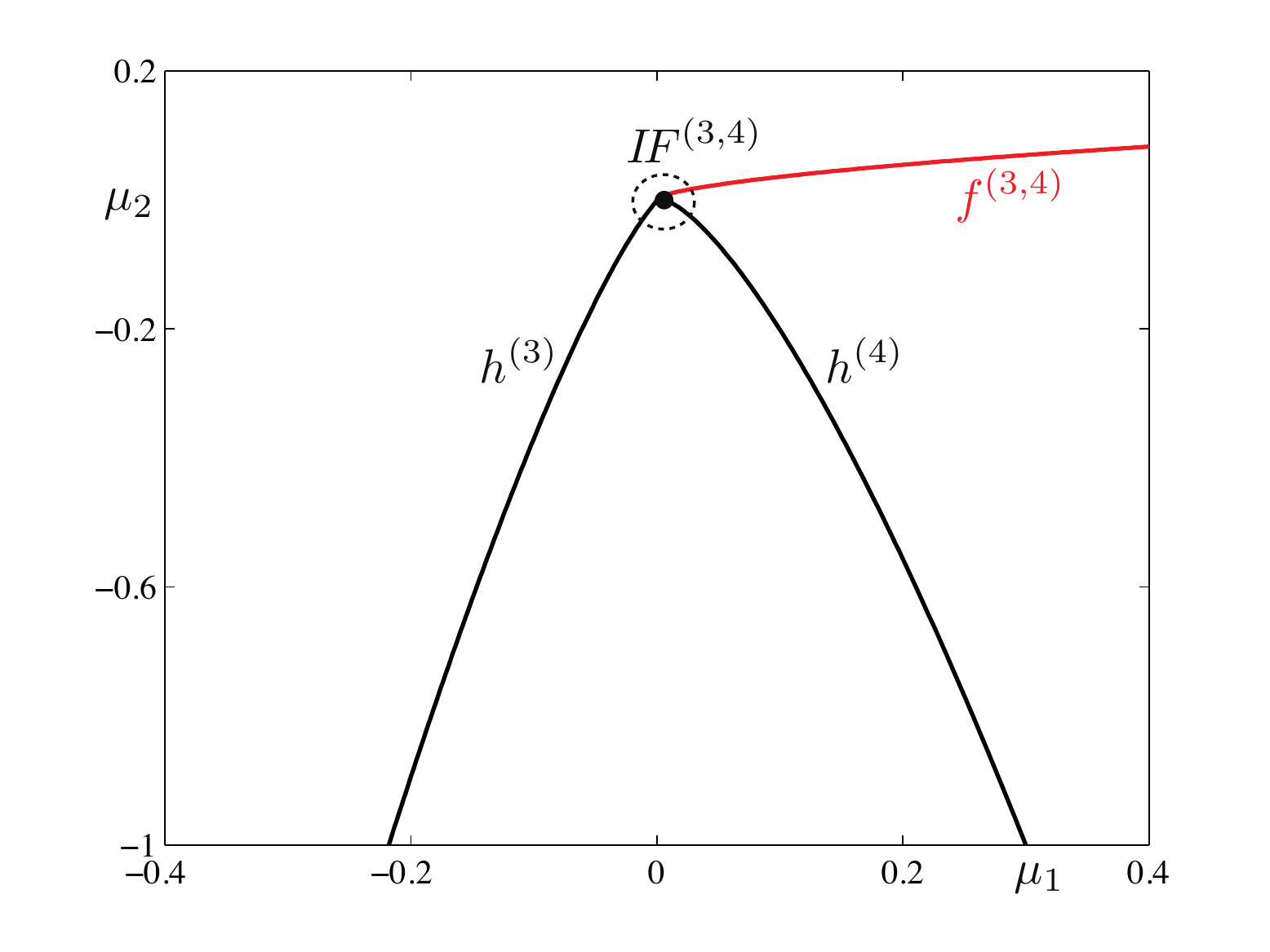}}
\caption{\label{fig:twoparbifmap}Curve of fold and homoclinic bifurcation points in the $(\mu_1,\mu_2)$-plane, obtained from the return map that we derived above.}
\end{figure}

We look for $\mu_1\minusone$ and $\mu_2\minusone$ families of fixed point of the previous set of
equations with \textsc{Auto}~\cite{Auto07PMan}. We fix the signs and continue in $\mu_1$ and in
$\mu_2$ the solutions to the system
\begin{eqnarray*}
\mu_1+K_1\xi_1^{\Delta_1}          &=& \xi_3-\xi_1\\
\mu_2\xi_3+K_2\xi_1^{\Delta_2}\xi_3^2 &=& \beta\xi_3^3+\frac{\xi_1}{2\beta}.
\end{eqnarray*}
We do find saddle-node bifurcations of equilibria which we can continue in two parameter and obtain
a curve of fold points in $(\mu_1,\mu_2)$-plane. In order to compute the curve of homoclinic
bifurcations in this plane, we use the fact that, in the map $\Pi_2$, the homoclinic connection
corresponds to $$\eta_1=\mu_1,\;\eta_2=\mu_2.$$ In the preceding system of equations, this gives:
$$K_1\xi_1^{\Delta_1}=0,\;K_2\xi_1^{\Delta_2}\xi_3=0.$$ Hence, we obtain the equation for the
homoclinic curve
\begin{eqnarray*}
\mu_1 &=& \xi_3-\xi_1\\
\mu_2 &=& \beta\xi_3^2+\frac{\xi_1}{2\beta\xi_3}.
\end{eqnarray*}
We can then compare the computed curve of folds with the curve of homoclinic points and we present
the result in figure~\ref{fig:twoparbifmap}. We obtain a qualitative agreement with the similar
curves computed from the HR system. Indeed, the homoclinic curve is folded and, from the tip of
that curve, corresponding to the inclination flip bifurcation $I\!F^{(3,4)}$, emanates a curve of
fold bifurcation, which corresponds to $f^{(3,4)}$. Note that our numerics is not valid in the
vicinity of this tip (dashed circle in~figure~\ref{fig:twoparbifmap}); however, the trend of both
the homoclinic and the fold curves outside this small region seems to indicate that they indeed
meet at the tip.

For values of $I$ and $b$ corresponding to the numerical return map described above and to the
inclination flip bifurcation, we can compute the eigenvalue ratios
$\Delta_2=-\lambda_{\!uu}/\lambda_{\!s}$ and $\Delta_1=-\lambda_{\!u}/\lambda_{\!s}$ and check
where the point $(\Delta_2,\Delta_1)$ is located in the diagram of Fig. 4 (left)
in~\cite{HomburgKrauskopf}, where different unfoldings of the inclination flip bifurcation are
studied. It appears that the HR system for the parameter values mentioned above falls into the case
C of the classification derived in~\cite{HomburgKrauskopf}; therefore, horseshoe dynamics is
expected in the vicinity of the inclination flip point, which is consistent with the results
of~\cite{Terman:1991}.

\section{Discussion}
\label{sec:5}
This paper has revisited the well-known Hindmarsh-Rose neuron model from a global bifurcation
analysis standpoint. To this end, we used different tools, geometrical and numerical. We
extracted specific information by relying on the strengths of each method and depicted a global
bifurcation scenario by exploiting the tools redundancy to overcome specific weaknesses of each
method.

In particular, the analysis we have carried out has shown that numerical continuation based on
boundary-value problems can be extremely useful in slow-fast systems where pure simulation can run
into difficulties in unfolding bifurcations that are very close to one another. Note also that the
homoclinic bifurcations studied here do not themselves represent stable dynamical behaviour.
Nevertheless, their influence on the global bifurcation structure is profound.

On the other hand, the geometrical analysis provided details that were not possible to obtain
numerically. Moreover, the geometrical analysis points to a generic phenomenon. Namely a sharp fold
in the curve of homoclinic orbits in the parameter plane should usually be associated with an
inclination flip. Such a phenomenon for example was also observed as part of the unfolding of a
tangent period-to-equilibrium heteroclinic cycle \cite{Champneys:2009}.

The obtained bifurcation scenario is organised by various curves of homoclinic bifurcations and
their codimension-two degeneracies and explains the smooth spike-adding transition (where the
number of spikes in each burst is increased by one) typical of the HR model and of many other
neuron models. In some sense, the work presented here extends the work of A. Shilnikov who also
detected the presence of IF and OF bifurcations in the HR neuron model. He found a wealth of
complex dynamics, however, his bifurcation analysis was limited to the curve $h^{(1)}$. The results
here have shown that the key to understanding the origins of the spike adding behaviour is to
analyse the IFs and OFs occurring on the homoclinic curves $h^{(n)}$ for $n>1$.

The analysis reported in this paper is interesting not only for its intrinsic value in explaining
spike-adding in the HR neuron model, but also because similar bifurcation structures have been
found and analysed in other studies, such as those reported in \cite{Channell:2007,Mosekilde:2001,
Shilnikov:2005}. In particular, in \cite{Mosekilde:2001} the authors perform a
bifurcation analysis of a model of pancreatic $\beta$-cells, which show excitable features similar
to those of neurons, and find a global bifurcation structure that strikingly resembles what we
found for the HR model (see Fig. 4 of the cited paper). In
\cite{Channell:2007,Shilnikov:2005}, the authors present an analysis of a reduced model of leech
heart interneuron: also in this case, the period-adding mechanism is regulated by the presence of
homoclinic bifurcations and their degeneracies.

Obviously, a detailed bifurcation analysis of each model will show differences among models, but we
dare say that the {\em global} bifurcation structure, \ie the presence of homoclinic bifurcations
and the interplay of period doubling and fold of cycles bifurcations, remains unchanged and
constitutes a trademark of models of excitable cells that display the widespread period-adding
mechanism.

Even if the global bifurcation scenario is quite clear, aspects of the dynamics that remain to be
investigated include:
\begin{itemize}
\item Further explanation in a neighbourhood of the $IF^{(1)}$ point, in particular whether the
    $IF^{(n)}$ are really the same bifurcation point or not;
\item What happens to the primary homoclinic orbit before $IF^{(1)}$;
\item A rigorous analysis using slow-fast methods of the bifurcations caused by the homoclinic
    orbit passing through a neighbourhood of the fold of the critical manifold;
\item Detailed investigation of these phenomena in other models.
\end{itemize}

\vspace{0.2cm}
\noindent
{\bf Acknowledgement}:
The research of M.D. was supported by EPSRC under grant EP/E032249/1.

\end{document}